\documentclass[12pt,twoside]{amsart}
\def\ver{xm-2i.tex}

\usepackage{amssymb}
\usepackage{amsmath}
\usepackage{amscd}
\usepackage{calc}
\usepackage{tikz}
\usepackage{pifont}
\usepackage{color}
\usepackage{cite}

\usetikzlibrary{cd}

\usepackage{graphicx}
\def\qed{$\s$}
\def\Gproj{\operatorname{G-proj}}
\def\uGproj{\operatorname{\underline{G-proj}}}
\def\mod{\operatorname{mod}}
\def\Hom{\operatorname{Hom}}
\def\ind{\operatorname{ind}}
\def\mon{\operatorname{mon}}
\def\Gp{\operatorname{Gp}}
\def\umod{\operatorname{\underline{mod}}}
\def\Cok{\operatorname{Cok}}

\def\Ima{\operatorname{Im}}  % new added
\def\op{{\rm op}}
\def\Mimo{{\rm Mimo}}
\def\Mono{{\rm Mono}}
\def\from{\leftarrow}
\renewcommand\bar{\overline}

\numberwithin{equation}{section}
\newtheorem{thm}[equation]{\sc Theorem}
\newtheorem{lem}[equation]{\sc Lemma}
\newtheorem{cor}[equation]{\sc Corollary}
\newtheorem{obs}[equation]{\sc Observation}
\newtheorem{prop}[equation]{\sc Proposition}

\newtheoremstyle{notation}{3pt}{3pt}{}{}{\itshape}{:}{.5em}{\thmname{#1}}
\theoremstyle{notation}

\newtheorem{rem}{\it Remark}

\newtheorem{defin}{\it Definition}
\newtheorem{ex}{\it Example}

\definecolor{darkgreen}{rgb}{0,0.5,0}

\newcommand{\s}{\hfill \blacksquare}
\newcommand{\D}{\displaystyle}
\newenvironment{smallpmatrix}{\left(\begin{smallmatrix}}{\end{smallmatrix}\right)}

\input prepictex   \input pictex    \input postpictex
 \def\arr#1#2{\arrow <2mm> [0.25,0.75] from #1 to #2}
 \def\sq{\plot 0 0  1 0  1 1  0 1  0 0 /}
 \def\bul{$\sssize\bullet$}
 \def\trr{$\sssize\blacktriangleleft$}
 \def\trl{$\sssize\blacktriangleright$}
 \def\cir{\circulararc 360 degrees from 0 0.7 center at 0 0.4 }
\def\sssize{\scriptscriptstyle}

\def\circled#1{\beginpicture\setcoordinatesystem units <1mm,1mm>  point at 0 -1.5
  \multiput {} at -2 -2  2 2 /
  \put {$\bigcirc$} at 0 0
  \put {$\sssize\sf#1$} at 0 0 \endpicture}

\def\lac{\multiput{} at 0 -1  0 4 /
  \plot 0 -1  0 3 /
  \multiput{\sq} at -1 0  -1 1  -1 2 /
}
\def\lbb{\multiput{} at 0 -1  0 4 /
  \plot 0 -1  0 3 /
  \multiput{\sq} at -1 .5  -1 1.5  /
}
\def\lbd{\multiput{} at 0 -1  0 4 /
  \plot 0 -1  0 3 /
  \multiput{\sq} at -1 0  -1 1  -1 2  0 0 /
  \multiput{\bul} at -.5 0  .5 0 /
  \plot -.5 0  .5 0 /
}
\def\lbf{\multiput{} at 0 -1  0 4 /
  \plot 0 -1  0 3 /
  \multiput{\sq} at 0 1 /
}
\def\lbh{\multiput{} at 0 -1  0 4 /
  \plot 0 -1  0 3 /
  \multiput{\sq} at -1 0  -1 1  0 0  0 1  0 2 /
  \multiput{\bul} at -.5 1  .5 1 /
  \plot -.5 1  .5 1 /
}
\def\lca{\multiput{} at 0 -1  0 4 /
  \plot 0 -1  0 3 /
  \multiput{\sq} at -2 1  -1 0  -1 1  -1 2  0 0  0 1  0 2 /
  \multiput{\bul} at -1.5 1  -.5 1  .5 1 /
  \plot -1.5 1  .5 1 /
}
\def\lcc{\multiput{} at 0 -1  0 4 /
  \plot 0 -1  0 3 /
  \multiput{\sq} at -1 .5  -1 1.5   0 .5 /
  \multiput{\bul} at -.5 .5  .5 .5 /
  \plot -.5 .5  .5 .5 /
}
\def\lce{\multiput{} at 0 -1  0 4 /
  \plot 0 -1  0 3 /
  \multiput{\sq} at -1 0  -1 1  -1 2  0 0  0 1 /
  \multiput{\bul} at -.5 0  .5 0 /
  \plot -.5 0  .5 0 /
}
\def\lcg{\multiput{} at 0 -1  0 4 /
  \plot 0 -1  0 3 /
  \multiput{\sq} at -1 0  -1 1    0 0  0 1  0 2  1 1 /
  \multiput{\bul} at -.5 1  .5 1  1.5 1 /
  \plot -.5 1  1.5 1 /
}
\def\ldd{\multiput{} at 0 -1  0 4 /
  \plot 0 -1  0 3 /
  \multiput{\sq} at -1 0  -1 1  -1 2  0 0  0 1  0 2 /
  \multiput{\bul} at -.5 2  .5 2 /
  \plot -.5 2  .5 2 /
}
\def\lea{\multiput{} at 0 -1  0 4 /
  \plot 0 -1  0 3 /
  \multiput{\sq} at  -1 1   0 1 /
  \multiput{\bul} at -.5 1  .5 1 /
  \plot -.5 1  .5 1 /
}
\def\lec{\multiput{} at 0 -1  0 4 /
  \plot 0 -1  0 3 /
  \multiput{\sq} at -1 0  -1 1  -1 2  0 0  0 1  0 2 /
  \multiput{\bul} at -.5 1  .5 1 /
  \plot -.5 1  .5 1 /
}
\def\lee{\multiput{} at 0 -1  0 4 /
  \plot 0 -1  0 3 /
  \multiput{\sq} at -1 .5  -1 1.5   0 .5  0 1.5 /
  \multiput{\bul} at -.5 1.5  .5 1.5 /
  \plot -.5 1.5  .5 1.5 /
}
\def\leg{\multiput{} at 0 -1  0 4 /
  \plot 0 -1  0 3 /
  \multiput{\sq} at -1 0  -1 1  -1 2  0 0  0 1  0 2 /
  \multiput{\bul} at -.5 0  .5 0 /
  \plot -.5 0  .5 0 /
}
\def\lfb{\multiput{} at 0 -1  0 4 /
  \plot 0 -1  0 3 /
  \multiput{\sq} at -2 1  -1 0  -1 1  -1 2  0 0  0 1  0 2  1 1  /
  \multiput{\bul} at -1.7 1.2  -.7 1.2  .3 1.2  -.3 .8  .7 .8  1.7 .8 /
  \plot -1.7 1.2  .3 1.2 /
  \plot -.3 .8  1.7 .8 /
}
\def\lfd{\multiput{} at 0 -1  0 4 /
  \plot 0 -1  0 3 /
  \multiput{\sq} at -1 .5  -1 1.5  0 .5  0 1.5 /
  \multiput{\bul} at -.5 .5  .5 .5 /
  \plot -.5 .5  .5 .5 /
}
\def\lff{\multiput{} at 0 -1  0 4 /
  \plot 0 -1  0 3 /
  \multiput{\sq} at -2 0  -2 1  -2 2  -1 0  -1 1  0 0  0 1  1 0  1 1  1 2 /
  \multiput{\bul} at -.5 1  .5 1  -1.5 0  -.5 0  1.5 0 /
  \plot -.5 1  .5 1 /
  \plot -1.5 0  1.5 0 /
}
\def\lfh{\multiput{} at 0 -1  0 4 /
  \plot 0 -1  0 3 /
  \multiput{\sq} at -2 1  -1 0  -1 1  -1 2  0 0  0 1  0 2  1 1  /
  \multiput{\bul} at -1.5 1  -.5 1  .5 1  1.5 1 /
  \plot -1.5 1  1.5 1 /
}
\def\lga{\multiput{} at 0 -1  0 4 /
  \plot 0 -1  0 3 /
  \multiput{\sq} at -1 0  -1 1  -1 2  0 0  0 1  0 2  1 1 /
  \multiput{\bul} at -.5 1  .5 1  1.5 1 /
  \plot -.5 1  1.5 1 /
}
\def\lgc{\multiput{} at 0 -1  0 4 /
  \plot 0 -1  0 3 /
  \multiput{\sq} at -1 .5  0 .5  0 1.5   /
  \multiput{\bul} at -.5 .5  .5 .5 /
  \plot -.5 .5  .5 .5 /
}
\def\lge{\multiput{} at 0 -1  0 4 /
  \plot 0 -1  0 3 /
  \multiput{\sq} at -1 0  -1 1  0 0  0 1  0 2 /
  \multiput{\bul} at -.5 0  .5 0 /
  \plot -.5 0  .5 0 /
}
\def\lgg{\multiput{} at 0 -1  0 4 /
  \plot 0 -1  0 3 /
  \multiput{\sq} at -2 1  -1 0  -1 1  -1 2  0 0  0 1 /
  \multiput{\bul} at -1.5 1  -.5 1  .5 1 /
  \plot -1.5 1  .5 1 /
}
\def\lhb{\multiput{} at 0 -1  0 4 /
  \plot 0 -1  0 3 /
  \multiput{\sq} at  0 .5  0 1.5 /
}
\def\lhd{\multiput{} at 0 -1  0 4 /
  \plot 0 -1  0 3 /
  \multiput{\sq} at -1 0  0 0  0 1  0 2 /
  \multiput{\bul} at -.5 0  .5 0 /
  \plot -.5 0  .5 0 /
}
\def\lhf{\multiput{} at 0 -1  0 4 /
  \plot 0 -1  0 3 /
  \multiput{\sq} at -1 1   /
}
\def\lhh{\multiput{} at 0 -1  0 4 /
  \plot 0 -1  0 3 /
  \multiput{\sq} at -1 0  -1 1  -1 2  0 0  0 1 /
  \multiput{\bul} at -.5 1  .5 1 /
  \plot -.5 1  .5 1 /
}
\def\lic{\multiput{} at 0 -1  0 4 /
  \plot 0 -1  0 3 /
  \multiput{\sq} at   0 0  0 1  0 2 /
}

\def\arqkLambdaZeta{\beginpicture\setcoordinatesystem units <1.5cm,1.8cm>
	\put{\scale{\put{\lac}  at 0 0  }} at  2 6
	\put{\scale{\put{\lbb}  at 0 0  }} at  1 5
        \put{\scale{\put{\circled5}  at 0 .5  }} at  1.22 5.0
	\put{\scale{\put{\lbd}  at 0 0  }} at  3 5
	\put{\scale{\put{\lbf}  at 0 0  }} at  5 5
        \put{\scale{\put{\circled{11}}  at 0 .5  }} at  5.3 5.1
        \put{\scale{\put{\lbh}  at 0 0  }} at  7 5
	\multiput{\scale{\put{\lca}  at 0 0  }} at  0 4  8 4 /
	\put{\scale{\put{\lcc}  at 0 0  }} at  2 4
	\put{\scale{\put{\lce}  at 0 0  }} at  4 4
	\put{\scale{\put{\lcg}  at 0 0  }} at  6 4
	\put{\scale{\put{\ldd}  at 0 0  }} at  3 3.6
	\multiput{\scale{\put{\lea}  at 0 0  }} at  0 3.3  8 3.3 /
	\put{\scale{\put{\lec}  at 0 0  }} at  2 3.3
	\put{\scale{\put{\lee}  at 0 0  }} at  4 3.3
	\put{\scale{\put{\leg}  at 0 0  }} at  6 3.3
	\put{\scale{\put{\lfb}  at 0 0  }} at  1 3
	\put{\scale{\put{\lfd}  at 0 0  }} at  3 3
	\put{\scale{\put{\lff}  at 0 0  }} at  5 3
	\put{\scale{\put{\lfh}  at 0 0  }} at  7 3
	\multiput{\scale{\put{\lga}  at 0 0  }} at  0 2  8 2 /
	\put{\scale{\put{\lgc}  at 0 0  }} at  2 2
	\put{\scale{\put{\lge}  at 0 0  }} at  4 2
	\put{\scale{\put{\lgg}  at 0 0  }} at  6 2
	\put{\scale{\put{\lhb}  at 0 0  }} at  1 1
	\put{\scale{\put{\lhd}  at 0 0  }} at  3 1
        \put{\scale{\put{\circled8}  at 0 .5  }} at  3.4 .6
	\put{\scale{\put{\lhf}  at 0 0  }} at  5 1
	\put{\scale{\put{\lhh}  at 0 0  }} at  7 1
        \put{\scale{\put{\circled 2}  at 0 .5  }} at  7.4 .6
	\put{\scale{\put{\lic}  at 0 0  }} at  2 0

        \multiput{ \arr{.3 .35} {.6 .65} } at 2 .3  0 2.3  2 2.3  4 2.3  6 2.3 0 4.3  2 4.3  4 4.3  6 4.3
        1 1.3  3 1.3  5 1.3   7 1.3  1 3.3  3 3.3  5 3.3  7 3.3  1 5.3 /
        \multiput{ \arr{.3 .65} {.6 .35} } at  1 .3  0 1.3  2 1.3  4 1.3  6 1.3
        1 2.3  3 2.3  5 2.3  7 2.3  0 3.3  2 3.3  4 3.3  6 3.3  1 4.3  3 4.3  5 4.3  7 4.3  2 5.3 /
        \multiput{ \arr{.35 .12} {.55 .18} } at 1 3.1  3 3.1  5 3.1  7 3.1  2 3.4 /
        \multiput{ \arr{.35 .18} {.55 .12} } at 0 3.1  2 3.1  4 3.1  6 3.1  3 3.4 /

        \setdashes<3pt>
        \multiput{ \plot 0 0  0 1 /  \plot 0 2  0 2.3 /  \plot 0 4  0 5 / } at -.05 3  7.95 3 /

        \setdots<2pt>
        \plot 0 1  .5 1 /
        \plot 3.5 1  4.5 1 /
        \plot 5.5 1  6.5 1 /
        \plot 7.5 1  8 1 /
        \plot .5 3.3  1.5 3.3 /
        \plot 4.5 3.3  5.5 3.3 /
        \plot 6.5 3.3  7.5 3.3 /
        \plot 0 5  .5 5 /
        \plot 3.5 5  4.5 5 /
        \plot 5.5 5  6.5 5 /
        \plot 7.5 5  8 5 /
\endpicture}

\def\arqkLambdaEps{\beginpicture\setcoordinatesystem units <1.5cm,1.3cm>
	\multiput{\scale{\put{\lhb} at -.5 1}} at  0 4  4 0 /
	\multiput{\scale{\put{\lbb}  at .5 1}} at  0 0  4 4 /
        \put{\scale{\put{\lgc}  at 0 1  }} at  1 3
	\put{\scale{\put{\lhf}  at 0.5 1  }} at  3 3
        \multiput{\scale{\put{\lea}  at 0 1  }} at  0 2  4 2 /
	\put{\scale{\put{\lfd}  at 0 1  }} at  2 2
	\put{\scale{\put{\lcc}  at 0 1 }} at  1 1
	\put{\scale{\put{\lbf}  at -0.5 1 }} at  3 1
	\put{\scale{\put{\lee}  at 0 0.2 }} at  3 2.3
        \multiput{ \arr{.3 .35} {.6 .65} } at 0 0  1 1  3 1  0 2  2 2
           3 3 /
        \multiput{ \arr{.3 .65} {.6 .35} } at  3 0  0 1  2 1  1 2  3 2
           0 3 /
        \put{ \arr{.3 .09} {.6 .18} } at 2 1.8
        \put{ \arr{.3 .18} {.6 .09} } at 3 1.8

        \setdashes<3pt>
        \multiput{ \plot 0 -1  0 -.6  / \plot 0 0.7  0 1.4 /  \plot 0 2.7  0 3.4 /
          \plot 0 4.7  0 5 / } at -.05 1.7  3.95 1.7 /

        \setdots<2pt>
        \plot 1.5 .7  2.5 .7 /
        \plot 1.5 2.7  2.5 2.7 /
\endpicture}

\def\vag{\multiput{} at 0 -1  0 4 /
  \multiput{\sq} at 0 0  0 1  0 2  /
  \put{\trr} at .5 2
}
\def\vbb{\multiput{} at 0 -1  0 4 /
  \multiput{\sq} at 0 0  0 1  0 2 /
  \put{\trl} at .5 0
}
\def\vbd{\multiput{} at 0 -1  0 4 /
  \multiput{\sq} at 0 1 /
  \put{\trl} at .5 1
}
\def\vbf{\multiput{} at 0 -1  0 4 /
  \multiput{\sq} at 0 0  0 1  0 2 /
  \put{\trr} at .5 1
}
\def\vbh{\multiput{} at 0 -1  0 4 /
  \multiput{\sq} at 0 0  0 1 /
  \put{\trr} at .5 1
}
\def\vca{\multiput{} at 0 -1  0 4 /
  \multiput{\sq} at 0 0  0 1  0 2  1 0  1 1 /
  \multiput{\trl} at .5 0  1.5 0 /
  \put{\trr} at 1.5 1
  \plot .5 0  1.5 0 /
}
\def\vcc{\multiput{} at 0 -1  0 4 /
  \multiput{\sq} at 0 0  0 1  0 2  1 1 /
  \multiput{\trl} at .5 1  1.5 1 /
  \plot .5 1  1.5 1 /
}
\def\vce{\multiput{} at 0 -1  0 4 /
  \multiput{\sq} at 0 0  0 1  0 2  1 1 /
  \multiput{\trr} at .6 .8  1.6 .8 /
  \put{\trl} at 1.4 1.2
  \plot .6 .8  1.6 .8 /
}
\def\vcg{\multiput{} at 0 -1  0 4 /
  \multiput{\sq} at 0 0  0 1 /
  \put{\trr} at .5 0
}
\def\vdf{\multiput{} at 0 -1  0 4 /
  \multiput{\sq} at 0 0  0 1  0 2 /
}
\def\vea{\multiput{} at 0 -1  0 4 /
  \multiput{\sq} at 0 1 /
}
\def\vec{\multiput{} at 0 -1  0 4 /
  \multiput{\sq} at 0 0  0 1  0 2  1 0  1 1 /
  \multiput{\trl} at .4 1.2  1.4 1.2 /
  \put{\trr} at 1.6 .8
  \plot .4 1.2  1.4 1.2 /
}
\def\vee{\multiput{} at 0 -1  0 4 /
  \multiput{\sq} at 0 0  0 1 /
}
\def\veg{\multiput{} at 0 -1  0 4 /
  \multiput{\sq} at 0 0  0 1  0 2  1 0 /
  \multiput{\trl} at .4 .2  1.4 .2 /
  \put{\trr} at 1.6 -.2
  \plot .4 .2  1.4 .2 /
}
\def\vfb{\multiput{} at 0 -1  0 4 /
  \multiput{\sq} at 0 0  0 1  0 2  1 1  2 0  2 1 /
  \multiput{\trl} at .4 1.2  2.4 1.2 /
  \multiput{\trr} at 1.6 .8  2.6 .8 /
  \plot .4 1.2  2.4 1.2 /
  \plot 1.6 .8  2.6 .8 /
}
\def\vfd{\multiput{} at 0 -1  0 4 /
  \multiput{\sq} at 0 0  0 1  0 2  1 1  2 0  2 1  2 2 /
  \multiput{\trl} at .4 1.2  1.4 1.2 /
  \multiput{\trr} at 1.6 .8  2.6 .8 /
  \plot .4 1.2  1.4 1.2 /
  \plot 1.6 .8  2.6 .8 /
}
\def\vff{\multiput{} at 0 -1  0 4 /
  \multiput{\sq} at 0 0  0 1  1 0 /
  \multiput{\trl} at .4  .2  1.4 .2 /
  \put{\trr} at 1.6 -.2
  \plot .4 .2  1.4 .2 /
}
\def\vfh{\multiput{} at 0 -1  0 4 /
  \multiput{\sq} at 0 0  0 1  0 2  1 0  1 1 /
  \multiput{\trl} at .4 .2  1.4 .2 /
  \put{\trr} at 1.6 -.2
  \plot .4 .2  1.4 .2 /
}
\def\vga{\multiput{} at 0 -1  0 4 /
  \multiput{\sq} at 0 0  0 1  0 2  1 0  1 1 /
  \multiput{\trl} at .5 1  1.5 1 /
  \put{\trr} at 1.5 0
  \plot .5 1  1.5 1 /
}
\def\vgc{\multiput{} at 0 -1  0 4 /
  \multiput{\sq} at 0 0  0 1  0 2  1 1 /
  \multiput{\trr} at .5 1  1.5 1 /
  \plot .5 1  1.5 1 /
}
\def\vge{\multiput{} at 0 -1  0 4 /
  \multiput{\sq} at 0 0  0 1  0 2  1 1 /
  \multiput{\trl} at .4 1.2  1.4 1.2 /
  \put{\trr} at 1.6 .8
  \plot .4 1.2  1.4 1.2 /
}
\def\vgg{\multiput{} at 0 -1  0 4 /
  \multiput{\sq} at 0 0  0 1 /
  \put{\trl} at .5 0
}
\def\vhb{\multiput{} at 0 -1  0 4 /
  \multiput{\sq} at 0 0  0 1  0 2 /
  \put{\trr} at .5 0
}
\def\vhd{\multiput{} at 0 -1  0 4 /
  \multiput{\sq} at 0 1 /
  \put{\trr} at .5 1
}
\def\vhf{\multiput{} at 0 -1  0 4 /
  \multiput{\sq} at 0 0  0 1  0 2 /
  \put{\trl} at .5 1
}
\def\vhh{\multiput{} at 0 -1  0 4 /
  \multiput{\sq} at 0 0  0 1 /
  \put{\trl} at .5 1
}
\def\vig{\multiput{} at 0 -1  0 4 /
  \multiput{\sq} at 0 0  0 1  0 2 /
  \put{\trl} at .5 2
}

\def\arqkVZeta{\beginpicture\setcoordinatesystem units <1.5cm,1.8cm>
	\put{\scale{\put{\vag}  at 1 0  }} at  6 6
	\put{\scale{\put{\vbb}  at 1 0  }} at  1 5
        \put{\scale{\put{\circled{7}}  at 0 .5  }} at  1.6 5.1
	\put{\scale{\put{\vbd}  at 1 0  }} at  3 5
	\put{\scale{\put{\vbf}  at 1 0  }} at  5 5
        \put{\scale{\put{\circled{1}}  at 0 .5  }} at  5.6 5
        \put{\scale{\put{\vbh}  at 1 0  }} at  7 5
	\multiput{\scale{\put{\vca}  at .5 0  }} at  0 4  8 4 /
	\put{\scale{\put{\vcc}  at .5 0  }} at  2 4
	\put{\scale{\put{\vce}  at .5 0  }} at  4 4
	\put{\scale{\put{\vcg}  at 1 .5  }} at  6 4
	\put{\scale{\put{\vdf}  at 1 -.5  }} at  5 3.6
	\multiput{\scale{\put{\vea}  at 1 0  }} at  0 3.3  8 3.3 /
	\put{\scale{\put{\vec}  at .5 0  }} at  2 3.3
	\put{\scale{\put{\vee}  at 1 .5  }} at  4 3.3
	\put{\scale{\put{\veg}  at .5 0  }} at  6 3.3
	\put{\scale{\put{\vfb}  at 0 0  }} at  1 3
	\put{\scale{\put{\vfd}  at 0 0  }} at  3 3
	\put{\scale{\put{\vff}  at 0 .5  }} at  5 3
	\put{\scale{\put{\vfh}  at 0 0  }} at  7 3
	\multiput{\scale{\put{\vga}  at .5 0  }} at  0 2  8 2 /
	\put{\scale{\put{\vgc}  at .5 0  }} at  2 2
	\put{\scale{\put{\vge}  at .5 0  }} at  4 2
	\put{\scale{\put{\vgg}  at 1 .5  }} at  6 2
	\put{\scale{\put{\vhb}  at 1 0  }} at  1 1
	\put{\scale{\put{\vhd}  at 1 0  }} at  3 1
        \put{\scale{\put{\circled{10}}  at 0 .5  }} at  3.6 .6
	\put{\scale{\put{\vhf}  at 1 0  }} at  5 1
	\put{\scale{\put{\vhh}  at 1 0  }} at  7 1
        \put{\scale{\put{\circled{4}}  at 0 .5  }} at  7.6 .6
	\put{\scale{\put{\vig}  at 1 0  }} at  6 0

        \multiput{ \arr{.3 .35} {.6 .65} } at 6 .3  0 2.3  2 2.3  4 2.3  6 2.3 0 4.3  2 4.3  4 4.3  6 4.3
        1 1.3  3 1.3  5 1.3   7 1.3  1 3.3  3 3.3  5 3.3  7 3.3  5 5.3 /
        \multiput{ \arr{.3 .65} {.6 .35} } at  5 .3  0 1.3  2 1.3  4 1.3  6 1.3
        1 2.3  3 2.3  5 2.3  7 2.3  0 3.3  2 3.3  4 3.3  6 3.3  1 4.3  3 4.3  5 4.3  7 4.3  6 5.3 /
        \multiput{ \arr{.35 .12} {.55 .18} } at 1 3.1  3 3.1  5 3.1  7 3.1  4 3.4 /
        \multiput{ \arr{.35 .18} {.55 .12} } at 0 3.1  2 3.1  4 3.1  6 3.1  5 3.4 /

        \setdashes<3pt>
        \multiput{ \plot 0 0  0 1 /  \plot 0 2  0 2.3 /  \plot 0 4  0 5 / } at -.05 3  7.95 3 /

        \setdots<2pt>
        \plot 0 1  .5 1 /
        \plot 1.5 1  2.5 1 /
        \plot 3.5 1  4.5 1 /
        \plot 7.5 1  8 1 /
        \plot .5 3.3  1.5 3.3 /
        \plot 2.5 3.3  3.5 3.3 /
        \plot 6.5 3.3  7.5 3.3 /
        \plot 0 5  .5 5 /
        \plot 1.5 5  2.5 5 /
        \plot 3.5 5  4.5 5 /
        \plot 7.5 5  8 5 /
\endpicture}

\def\arqkVEps{\beginpicture\setcoordinatesystem units <1.5cm,1.3cm>
	\multiput{\scale{\put{\vhh} at -.5 1.5}} at  0 4  4 0 /
	\multiput{\scale{\put{\vbh}  at -.5 1.5}} at  0 0  4 4 /
        \put{\scale{\put{\vbd}  at -.5 1  }} at  1 3
	\put{\scale{\put{\vcg}  at -.5 1.5  }} at  3 3
        \multiput{\scale{\put{\vea}  at -.5 1  }} at  0 2  4 2 /
	\put{\scale{\put{\vff}  at -.5 1.5  }} at  2 2
	\put{\scale{\put{\vhd}  at -.5 1 }} at  1 1
	\put{\scale{\put{\vgg}  at -.5 1.5 }} at  3 1
	\put{\scale{\put{\vee}  at -.5 1.5 }} at  1 2.3
        \multiput{ \arr{.3 .35} {.6 .65} } at 0 0  1 1  3 1  0 2  2 2
           3 3 /
        \multiput{ \arr{.3 .65} {.6 .35} } at  3 0  0 1  2 1  1 2  3 2
           0 3 /
        \put{ \arr{.3 .09} {.6 .18} } at 0 1.8
        \put{ \arr{.3 .18} {.6 .09} } at 1 1.8

        \setdashes<3pt>
        \multiput{ \plot 0 -1  0 -.5  / \plot 0 0.5  0 1.5 /  \plot 0 2.5  0 3.5 /
          \plot 0 4.5  0 5 / } at -.05 1.7  3.95 1.7 /

        \setdots<2pt>
        \plot 1.5 .7  2.5 .7 /
        \plot 1.5 2.7  2.5 2.7 /
\endpicture}

\def\uae{\multiput{} at 0 0  0 4 /
  \multiput{\sq} at 0 0  0 1  0 2 /
  \put{\cir} at .5 2
}
\def\ubb{\multiput{} at 0 0  0 4 /
  \put{\sq} at 0 1
}
\def\ubd{\multiput{} at 0 0  0 4 /
  \multiput{\sq} at 0 0  0 1  0 2 /
  \put{\cir} at .5 1
}
\def\ubf{\multiput{} at 0 0  0 4 /
  \multiput{\sq} at 0 0  0 1  0 2 /
  \put{\cir} at .5 2
  \put{\bul} at .5 0
}
\def\ubh{\multiput{} at 0 0  0 4 /
  \put{\sq} at 0 1
  \put{\cir} at .5 1
  \put{\bul} at .5 1
}
\def\uca{\multiput{} at 0 0  0 4 /
  \multiput{\sq} at 0 0  0 1 /
  \put{\cir} at .5 0
  \put{\bul} at .5 0
}
\def\ucc{\multiput{} at 0 0  0 4 /
  \multiput{\sq} at 0 0  0 1  0 2  1 1 /
  \multiput{\cir} at .5 1  1.5 1 /
  \plot .5 1.3  1.5 1.3 /
}
\def\uce{\multiput{} at 0 0  0 4 /
  \multiput{\sq} at 0 0  0 1  0 2 /
  \put{\cir} at .5 1
  \put{\bul} at .5 0
}
\def\ucg{\multiput{} at 0 0  0 4 /
  \multiput{\sq} at 0 0  0 1  0 2  1 1 /
  \multiput{\cir} at .5 2  1.5 1 /
  \multiput{\bul} at .5 1  1.5 1 /
  \plot .5 1  1.5 1 /
}
\def\uda{\multiput{} at 0 0  0 4 /
  \multiput{\sq} at 0 0  0 1 /
  \put{\cir} at .5 1
}
\def\udc{\multiput{} at 0 0  0 4 /
  \multiput{\sq} at 0 0  0 1  0 2 /
  \put{\cir} at .5 0
  \put{\bul} at .5 0
}
\def\ude{\multiput{} at 0 0  0 4 /
  \put{\sq} at 0 0
  \put{\cir} at .5 0
}
\def\udg{\multiput{} at 0 0  0 4 /
  \multiput{\sq} at 0 0  0 1  0 2  /
  \put{\cir} at .5 1
  \put{\bul} at .5 1
}
\def\ueb{\multiput{} at 0 0  0 4 /
  \multiput{\sq} at 0 0  0 1  1 0  1 1  1 2 /
  \multiput{\cir} at .5 1  1.5 0 /
  \multiput{\bul} at .5 0  1.5 0 /
  \plot .5 0  1.5 0 /
}
\def\ued{\multiput{} at 0 0  0 4 /
  \multiput{\sq} at 0 0  0 1  0 2  1 1 /
  \multiput{\cir} at .5 1  1.5 1 /
  \put{\bul} at .5 0
  \plot .5 1.3  1.5 1.3 /
}
\def\uef{\multiput{} at 0 0  0 4 /
  \multiput{\sq} at 0 0  0 1  0 2  1 1 /
  \multiput{\cir} at .5 1  1.5 1 /
  \multiput{\bul} at .5 1  1.5 1 /
  \plot .5 1  1.5 1 /
}
\def\ueh{\multiput{} at 0 0  0 4 /
  \multiput{\sq} at 0 0  0 1  0 2  1 1  1 2 /
  \multiput{\cir} at .5 2  1.5 1 /
  \multiput{\bul} at .5 1  1.5 1 /
  \plot .5 1  1.5 1 /
}
\def\ufa{\multiput{} at 0 0  0 4 /
  \multiput{\sq} at 0 0  0 1  0 2  1 1  1 2  1 3 /
  \multiput{\cir} at .5 2  1.5 1 /
  \multiput{\bul} at .5 1  1.5 1 /
  \plot .5 1  1.5 1 /
}
\def\ufc{\multiput{} at 0 0  0 4 /
  \multiput{\sq} at 0 0  0 1 /
  \put{\cir} at .5 1
  \put{\bul} at .5 0
}
\def\ufe{\multiput{} at 0 0  0 4 /
  \multiput{\sq} at 0 0  0 1  0 2  1 1 /
  \multiput{\cir} at .5 1  1.5 1 /
  \multiput{\bul} at .5 1  1.5 1 /
  \plot .5 1  1.5 1 /
  \plot .5 1.3  1.5 1.3 /
}
\def\ufg{\multiput{} at 0 0  0 4 /
  \multiput{\sq} at 0 0  0 1 /
  \put{\cir} at .5 0
}
\def\ugb{\multiput{} at 0 0  0 4 /
  \multiput{\sq} at 0 0  0 1  0 2 /
  \put{\cir} at .5 2
  \put{\bul} at .5 1
}
\def\ugd{\multiput{} at 0 0  0 4 /
  \multiput{\sq} at 0 0  0 1 /
  \put{\cir} at .5 1
  \put{\bul} at .5 1
}
\def\ugf{\multiput{} at 0 0  0 4 /
  \multiput{\sq} at 0 0  0 1 /
}
\def\ugh{\multiput{} at 0 0  0 4 /
  \multiput{\sq} at 0 0  0 1  0 2 /
  \put{\cir} at .5 0
}
\def\uhc{\multiput{} at 0 0  0 4 /
  \multiput{\sq} at 0 0  0 1  0 2 /
  \put{\cir} at .5 2
  \put{\bul} at .5 2
}
\def\uhg{\multiput{} at 0 0  0 4 /
  \multiput{\sq} at 0 0  0 1  0 2 /
}

\def\arqkLZeta{\beginpicture\setcoordinatesystem units <1.5cm,1.8cm>
	\put{\scale{\put{\uae}  at 0 0  }} at  4 6
	\put{\scale{\put{\ubb}  at 0 .5  }} at  1 5
	\put{\scale{\put{\circled0,\circled{12}}  at 0 .5  }} at  1.5 5.05
	\put{\scale{\put{\ubd}  at 0 0  }} at  3 5
	\put{\scale{\put{\ubf}  at 0 0  }} at  5 5
	\put{\scale{\put{\circled6}  at 0 .5  }} at  5.3 5.1
        \put{\scale{\put{\ubh}  at 0 .5  }} at  7 5
	\multiput{\scale{\put{\uca}  at 0 .5  }} at  0 4  8 4 /
	\put{\scale{\put{\ucc}  at -1 0  }} at  2 4
	\put{\scale{\put{\uce}  at 0 .5  }} at  4 4
	\put{\scale{\put{\ucg}  at -1 0  }} at  6 4
	\multiput{\scale{\put{\uda}  at 0 .5  }} at  0 3.3  8 3.3 /
	\put{\scale{\put{\udc}  at 0 0  }} at  2 3.3
	\put{\scale{\put{\ude}  at 0 1  }} at  4 3.3
	\put{\scale{\put{\udg}  at 0 0  }} at  6 3.3
	\put{\scale{\put{\ueb}  at -1 0  }} at  1 3
	\put{\scale{\put{\ued}  at -1 0  }} at  3 3
	\put{\scale{\put{\uef}  at -1 0  }} at  5 3
	\put{\scale{\put{\ueh}  at -1 0  }} at  7 3
	\multiput{\scale{\put{\ufa}  at -1 -.5  }} at  0 2  8 2 /
	\put{\scale{\put{\ufc}  at 0 .5  }} at  2 2
	\put{\scale{\put{\ufe}  at -1 0  }} at  4 2
	\put{\scale{\put{\ufg}  at 0 .5  }} at  6 2
	\put{\scale{\put{\ugb}  at 0 0  }} at  1 1
	\put{\scale{\put{\ugd}  at 0 .5  }} at  3 1
        \put{\scale{\put{\circled{3}}  at 0 .5  }} at  3.4 .6
	\put{\scale{\put{\ugf}  at 0 .5  }} at  5 1
	\put{\scale{\put{\ugh}  at 0 0  }} at  7 1
        \put{\scale{\put{\circled{9}}  at 0 .5  }} at  7.4 .6
	\put{\scale{\put{\uhc}  at 0 0  }} at  2 0
	\put{\scale{\put{\uhg}  at 0 0  }} at  6 0

        \multiput{ \arr{.3 .35} {.6 .65} } at 2 .3  6 .3  0 2.3  2 2.3  4 2.3  6 2.3 0 4.3  2 4.3  4 4.3  6 4.3
        1 1.3  3 1.3  5 1.3   7 1.3  1 3.3  3 3.3  5 3.3  7 3.3  3 5.3 /
        \multiput{ \arr{.3 .65} {.6 .35} } at  1 .3  5 .3  0 1.3  2 1.3  4 1.3  6 1.3
        1 2.3  3 2.3  5 2.3  7 2.3  0 3.3  2 3.3  4 3.3  6 3.3  1 4.3  3 4.3  5 4.3  7 4.3  4 5.3 /
        \multiput{ \arr{.35 .12} {.55 .18} } at 1 3.1  3 3.1  5 3.1  7 3.1 /
        \multiput{ \arr{.35 .18} {.55 .12} } at 0 3.1  2 3.1  4 3.1  6 3.1 /

        \setdashes<3pt>
        \multiput{ \plot 0 0  0 1 /  \plot 0 2  0 2.3 /  \plot 0 4  0 5 / } at -.05 3  7.95 3 /

        \setdots<2pt>
        \plot 0 1  .5 1 /
        \plot 3.5 1  4.5 1 /
        \plot 7.5 1  8 1 /
        \plot .5 3.27  1.5 3.27 /
        \plot 2.5 3.27  3.5 3.27 /
        \plot 4.5 3.27  5.5 3.27 /
        \plot 6.5 3.27  7.5 3.27 /
        \plot 0 5  .5 5 /
        \plot 1.5 5  2.5 5 /
        \plot 5.5 5  6.5 5 /
        \plot 7.5 5  8 5 /
\endpicture}

\def\arqkLEps{\beginpicture\setcoordinatesystem units <1.5cm,1.3cm>
	\multiput{\scale{\put{\ugf} at -.5 1.5}} at  0 4  4 0 /
	\put{\scale{\put{\uda}  at -.5 1.5  }} at  2 4
	\multiput{\scale{\put{\ugd}  at -.5 1.5 }} at  0 0  4 4 /
        \put{\scale{\put{\ufg}  at -.5 1.5  }} at  1 3
	\put{\scale{\put{\ufc}  at -.5 1.5  }} at  3 3
        \multiput{\scale{\put{\ude}  at -.5 2  }} at  0 2  4 2 /
	\put{\scale{\put{\uca}  at -.5 1.5  }} at  2 2
	\put{\scale{\put{\ubh}  at -.5 1.5 }} at  1 1
	\put{\scale{\put{\ubb}  at -.5 1.5 }} at  3 1
        \multiput{ \arr{.3 .35} {.6 .65} } at 0 0  1 1  3 1  0 2  2 2
        1 3  3 3 /
        \multiput{ \arr{.3 .65} {.6 .35} } at  3 0  0 1  2 1  1 2  3 2
        0 3  2 3 /

        \setdashes<3pt>
        \multiput{ \plot 0 -1  0 -.5  / \plot 0 0.5  0 1.5 /  \plot 0 2.5  0 3.5 /
          \plot 0 4.5  0 5 / } at -.05 1.6  3.95 1.6 /

        \setdots<2pt>
        \plot 1.5 .6  2.5 .6 /
\endpicture}

\begin{document}

{[\footnotesize{\tt \ver,} \today]}\bigskip

\begin{center}
  {\large\bf
  A Reflection Equivalence\smallskip

  for Gorenstein-Projective Quiver Representations}

  \medskip
  Xiu-Hua Luo%
  \footnote{The first named author is supported by NSFC (No.\ 11771272) and
    Jiangsu Government Scholarship for Oversea Studies.}
  and Markus Schmidmeier%
  \footnote{This research is partially supported
    by a~Travel and Collaboration Grant from the Simons Foundation (Grant number
    245848 to the second named author).}

  \bigskip \parbox{10cm}{\footnotesize{\bf Abstract:}
    For $\Lambda$ a selfinjective algebra, and $Q$ a finite quiver without
    oriented cycles, the algebra $\Lambda Q$ is a Gorenstein algebra and
    the category $\Gproj\Lambda Q$ of Gorenstein-projective $\Lambda Q$-modules is
    a Frobenius category.  For a sink $v$ of $Q$, we define a functor
    $F(v) :
    \uGproj\Lambda Q\to \uGproj\Lambda Q(v)$  between the stable
      categories modulo projectives, where $Q(v)$ is obtained from $Q$
    by changing the direction of each arrow ending in $v$.
      The functor is given by an explicit construction
      on the level of objects and homomorphisms.
    Our main result states that $F(v)$ is an equivalence of categories.
    In the case where the underlying graph of $Q$ is a tree,
    we deduce that the stable category $\uGproj\Lambda Q$ does not depend on
    the orientation of $Q$.
    Moreover, if $Q$ is a quiver of type $\mathbb A_3$ and
    $\Lambda=k[T]/(T^n)$ the bounded polynomial algebra, we use the
    symmetry of the octahedron in the 
    octahedral axiom to verify that the composition of twelve reflections
    yields the identity on objects.
  }

\medskip \parbox{10cm}{\footnotesize{\bf  MSC 2020:}
  18G25 (relative homological algebra),
  16G20 (representations of quivers)
}

\medskip \parbox{10cm}{\footnotesize{\bf Key words:}
  Gorenstein-projective, separated monic, quiver representations,
  triangulated category, octahedral axiom

}
\end{center}

\section{Introduction}
Throughout this paper, let $k$ be a field,
$\Lambda$ be a finite-dimensional selfinjective $k$-algebra,
and $Q$ a finite quiver
  which is acyclic, that is, $Q$ has no oriented cycles.
The {\it path algebra} of the quiver $Q$ with coefficients in the algebra $\Lambda$
is denoted by $\Lambda Q=\Lambda\otimes_k kQ$.
Since $\Lambda$ is selfinjective, and $kQ$ hereditary, $\Lambda Q$
is a Gorenstein algebra \cite{AR2}; then the category $\Gproj\Lambda Q$ of Gorenstein-projective
modules is a Frobenius category.

\medskip
Let $Q$ be a quiver and let the vertex $v\in Q_0$ be a sink.  The {\it reflection} of $Q$ at $v$ is the quiver $Q(v)$
on the same set of vertices,
but with each arrow ending in $v$ replaced by its opposite.

  \medskip
  In this paper we use a pull-back construction on the level of
  objects and homomorphisms in the category $\Gproj\Lambda Q$
  to define the functor $F(v)$ in the following main result.

  \begin{thm}
    \label{theorem-intro1}

    For $ \Lambda $ a selfinjective algebra and $Q$ a finite acyclic quiver,
    the functor $$ F(v): \uGproj\Lambda Q\to \uGproj\Lambda Q(v) $$
    is an equivalence of categories.
  \end{thm}

  In \cite[ Corollaries~7.3--7.4 ]{ABHV},
  the authors proved that there is an equivalence of singularity categories $\mathcal{D}_{sg}(\Lambda Q)\cong \mathcal{D}_{sg}(\Lambda Q(v))$, where $\Lambda$ is a Gorenstein algebra. As a consequence, they
  obtain a different proof for the
  equivalence $\uGproj\Lambda Q\cong \uGproj\Lambda Q(v)$.
    Due to the approach via singularity categories, it is not immediately
    clear how the objects in $\uGproj\Lambda Q$ and in
    $\uGproj\Lambda Q(v)$ correspond to each other.
    By comparison, the functor $F(v)$ in our paper is given by an explicit construction. 
    
    \medskip
    The main tool in our approach is the characterization
    of Gorenstein-projective quiver representations by Luo and Zhang
    \cite{LZ1, LZ2}
    as the separated monic representations, see Section~\ref{section-sep-mon}
    below.
    In particular for a quiver of type $\mathbb A_n$ in linear operation,
    the Gorenstein-projective quiver representations are the chains
    of monomorphisms.  Thus, in the special case where
    $\Lambda=k[T]/(T^n)$ is the bounded polynomial ring,
    we are dealing with systems consisting of a nilpotent linear operator and
    a chain of $n-1$ invariant subspaces, see \cite{RS2} and \cite{KLM};
    such systems occur in applications, for example in control theory
    \cite{MS}.

    \smallskip
    Gorenstein-projective modules feature prominently in relative homological
    algebra where they assume the role played by projectives modules
    in homological algebra \cite{EJ1,EJ2}.  Ringel and Zhang exhibit a
    link to derived categories: If $\Lambda=k[\varepsilon]=k[T]/(T^2)$
    is the bounded polynomial ring, the authors show the equivalence
    $\uGproj \Lambda Q\simeq D^b(\mod kQ)/[1]$ in \cite{RZ}.

\bigskip
In the second part of this paper we deal with examples
where the quiver $Q$ has type $\mathbb A_3$
and $\Lambda$ is the bounded polynomial ring $k[T]/(T^n)$.
Up to quiver automorphism,
there are three orientations for the quiver:
$$Q:\; 1\to2\to3,\qquad Q(3):\; 1\to 2\leftarrow 3,\qquad
Q( 2,3 ):\;
1\leftarrow2\to 3.$$
So our result yields three categorical equivalences:
\begin{eqnarray*}
  F(3): & \uGproj \Lambda Q & \to\uGproj\Lambda Q(3),\\
  F(2,3): & \uGproj\Lambda Q(3) & \to \uGproj\Lambda Q(2,3),\\
  F( 1,2,3 ): & \uGproj\Lambda Q( 2,3 ) & \to \uGproj\Lambda Q.\end{eqnarray*}
The composition of three reflections yields a selfequivalence $F^3$ on $\uGproj\Lambda Q$.

\medskip

We show in Section~\ref{section-reflections}
that the composition of twelve reflections, $(F^3)^4$,
yields isomorphisms for objects:

\begin{thm}
  \label{theorem-intro2}

  Suppose $Q$ is a quiver of type $\mathbb A_3$ and $\Lambda=k[T]/(T^n)$.
  Let $M$ be a representation for $Q$ with coefficients in $\Lambda$ such that
  no branch $M_i$ has a nonzero projective direct summand.  Let $A$ be a right
  approximation for $M$ in $\Gproj\Lambda Q$.  Let $B$ be obtained from $A$ by
  applying twelve reflections. Then $$B\cong A\quad\text{in}\;\uGproj\Lambda Q.$$
\end{thm}

  The twelve reflections can be visualized on the octahedron in the octahedral axiom.
  Recall that of the eight faces of the octahedron, four correspond to exact triangles
  while the other four represent commutative triangles, each given by a pair of
  composable maps and their product.
  The four products form the cyclically oriented square in the octahedron which
  we consider the ``equator''.

  \medskip
  As an oriented graph, the octahedron has a symmetry group which is cyclic of order four
  (Observation~\ref{obs-symmetry}).
  A generator $\rho$ is given by rotating the octahedron by $90^\circ$ about the axis
  perpendicular to the equator, in the direction given by the four products,
  followed by a reflection on the plane of the equator.
  We show that $F(1,2,3)$ maps a composable pair to the pair in the octahedron
  obtained by applying $\rho$, up to the inverse of the suspension $[-1]$.
  Since $\rho$ has order four, any composition of twelve reflections yields
  the identity in the stable category, up to applications of the suspension.

  \medskip
  We conclude the paper with examples.

  \medskip
  First we exhibit the Auslander-Reiten quivers for the Gorenstein-pro\-jec\-tive modules over
  $\Lambda Q$, $\Lambda Q(3)$ and $\Lambda Q(2,3)$  where $\Lambda=k[T]/(T^2)$
  and $Q$ the above quiver of type $\mathbb A_3$ in linear orientation.
  While the three stable categories of Gorenstein-projective modules are equivalent,
  this result does not hold for the categories of all Gorenstein-projective modules,
  of all monic representations, or of all modules.

  \medskip
  Next we compute the orbit of a simple representation under reflections and illustrate
  it in the Auslander-Reiten quivers of
  $\Lambda Q$, $\Lambda Q(3)$ and $\Lambda Q(2,3)$
  for $\Lambda=k[T]/(T^3)$. As a byproduct we obtain the orbit in Figure~\ref{figure-one-orbit}.
  The icons will be explained in the text.
  The modules in the top row are representations of $\Lambda Q$, those in the second row
  of $\Lambda Q(2,3)$ and those in the bottom row of $\Lambda Q(3)$.

  \def\scale#1{\hbox{\beginpicture\setcoordinatesystem units <3mm,3mm>#1
    \multiput{} at -3 -4 3 3 / 
    \endpicture}}
\def\mugxm{\beginpicture\setcoordinatesystem units <.9cm,.7cm>
        \def\bul{$\scriptstyle\bullet$}
	\multiput{\scale{\put{\ued}  at -1 0  }} at 0 4  12 4 /
	\put{\scale{\put{\vfh}  at -1 0  }} at 1 0
	\put{\scale{\put{\lfb}  at 0 -.3  }} at 2 2
	\put{\scale{\put{\uef}  at -1 0  }} at 3 4
	\put{\scale{\put{\vfb}  at -1.5 0  }} at 4 0
	\put{\scale{\put{\lfd}  at 0 -.3  }} at 5 2
	\put{\scale{\put{\ueh}  at -1 0  }} at 6 4
	\put{\scale{\put{\vfd}  at -1.5 0  }} at 7 0
	\put{\scale{\put{\lff}  at 0 -.3 }} at 8 2
	\put{\scale{\put{\ueb}  at -1 0  }} at 9 4
	\put{\scale{\put{\vff}  at -1 0  }} at 10 0
	\put{\scale{\put{\lfh}  at 0 -.3  }} at 11 2
        \setquadratic
        \multiput {  \plot -.1 2.4  .05 .9  .3 .2 /
          \arr {.26 .26}  {.3 .2} } at 0 1.5  3 1.5  6 1.5  9 1.5 /
        \multiput { \plot 1.6 -.2  1.9 0  2.1 .8 /
          \arr {2.094 .76} {2.1 .8} } at 0 0  3 0  6 0  9 0 /
        \multiput { \plot 1.88 3  2.1 3.5   2.4 3.5 /
          \arr {2.36 3.52} {2.4 3.5} }  at  0 2.2  3 2.2  6 2.2  9 2.2 /
        \multiput {$\dots$} at -1 1.9  13 1.9 /
%        \put {{\tiny\tt X.L., M.S.}} at 11.8 -.8
%        \put{{\bf Gorenstein-projective}} at 2.5 -2.8
%        \put{{\bf Quiver Representations}} at 2.5 -3.5
%        \put{Top:\strut} at 7.5 -2.8
%        \put{Middle:\strut} at 9.5 -2.8
%        \put{Bottom:\strut} at 11.5 -2.8
%        \put{$\bullet\to\,\,\cirbox\,\,\to\Box$\strut} at 7.5 -3.5
%        \put{$\Box\leftarrow\bullet\to\Box$\strut} at 9.5 -3.5
%        \put{$\blacktriangleright\,\to\Box\leftarrow\,\blacktriangleleft$\strut}
%        at 11.5 -3.5
        \endpicture}

\begin{figure}[ht]
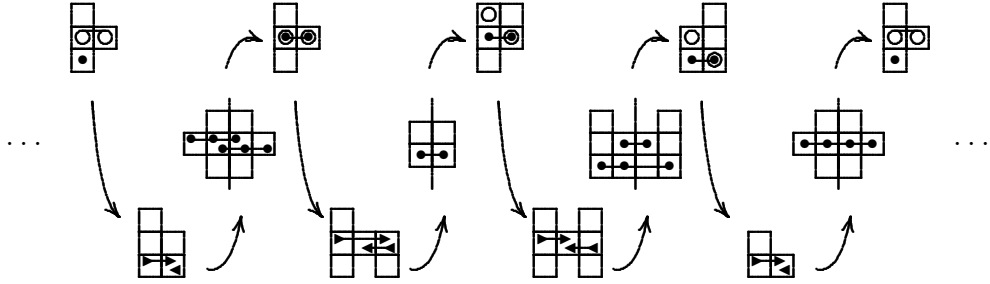

  \label{figure-one-orbit}
 
  \mugxm
  \caption{One orbit under the twelve reflections}
 
\end{figure}

\medskip
 
We summarize the contents of the paper. 

\medskip
In Section~\ref{section-monic} we study the interplay between the categories of
monic representations of $Q$ over $\Lambda$ and of
representations of $Q$ with objects and morphisms in the stable category
$\umod\Lambda$ modulo projectives.

\medskip The aim of Section~\ref{section-equivalence} is to present the
construction of objects and homomorphisms in $\Gproj\Lambda Q(v)$ from
those in $\Gproj\Lambda Q$.
Moreover, we sketch the proof of Theorem~\ref{theorem-intro1}.

\medskip
In Section~\ref{section-reflections} we revisit the octahedral axiom and use it
to trace the composition of twelve reflections on objects in
$\Gproj\Lambda Q$ where the quiver $Q$ has type $\mathbb A_3$,
showing Theorem~\ref{theorem-intro2}.

\medskip
The last Section~\ref{section-examples} gives examples which illustrate our
construction and visualize the composition of the twelve reflections.

%================================================================================
\section{Monic Representations}
%================================================================================
\label{section-monic}

%--------------------------------------------------------------------------------
\subsection{Notation}
%--------------------------------------------------------------------------------
\label{section-monic-notation}

We specify notation to define for $X\in\Gproj \Lambda Q $ and $v$ a sink in $Q$ the corresponding
module $X(v)\in \Gproj \Lambda Q(v) $.

\medskip
For a vertex $w\in Q_0$, let $w^-=\{s(\alpha)\;:\;\alpha\in Q_1\;\text{and}\;t(\alpha)= w \}$
be the multiset of predecessors of $w$.
Put
$$X_{w^-}=\bigoplus_{t(\alpha)=w}X_{s(\alpha)}.$$
We recall that a Gorenstein-projective module $X$ is
  separated monic (see definition  and results in Subsection~\ref{section-sep-mon}),
    so the map
    $$u_{X,w}=(X_\alpha)_{t(\alpha)=w}:X_{w^-}\to X_w$$
is a monomorphism, hence the short exact sequence
$$0\to X_{w^-}\stackrel {u_{X,w}}{\to} X_w\to X_{\bar w}\to 0$$
defines a $\Lambda$-module $X_{\bar w}$.  We will sometimes consider $X_{w^-}$ as a submodule of
$X_w$.  Regarding the map $u_{X,w}$, we sometimes omit one or both of the
subscripts.
\medskip

%--------------------------------------------------------------------------------
\subsection{Separated monic representations}
%--------------------------------------------------------------------------------
\label{section-sep-mon}

Given a finite  quiver $Q$ without oriented cycles, in their papers \cite{LZ1, LZ2},
the authors define separated monic representations of $Q$ over an algebra $\Lambda$ as follows:
\medskip

\begin{defin}  A representation $X = (X_i,X_{\alpha} )$ of $Q$ over $\Lambda$ is a
  {\it separated monic representation}, or a  {\it separated monic} $\Lambda Q$-module,
  if for each vertex $i\in Q_0$ the map
$$u_{X,i}: X_{i^-}\to   X_i $$
is a monomorphism. This means that $X$ satisfies the following two conditions:

\begin{itemize}
  \item[$(m1)$]   For each $\alpha\in Q_1$, $X_{\alpha}: X_{s(\alpha)}
\longrightarrow X_{t(\alpha)}$ is an injective homomorphism, and
\item[$(m2)$]  For each $i\in Q_0$,  the sum  $\sum\limits_{\begin
{smallmatrix}t(\alpha) =i \end{smallmatrix}}\Ima
X_\alpha  $ is a direct sum.
\end{itemize}
\end{defin}

If $\Lambda$ is a selfinjective algebra, then
the Gorenstein-projective representations of $Q$ over $\Lambda $
are exactly the separated monic representations \cite{LZ1, LZ2}.

%--------------------------------------------------------------------------------
\subsection{Minimal Monomorphisms}
%--------------------------------------------------------------------------------
\label{subsection-mimo}

  The algebra $\Lambda Q$ is Gorenstein, so every $\Lambda Q$-module has a right
  approximation in the category $\Gproj\Lambda Q$.

  \medskip
  Here we present a stepwise construction for such an approximation
  of a $\Lambda Q$-module $M$ in $\Gproj\Lambda Q$ which we call
  $\Mono(M)$ and define
  $\Mimo(M)$ to be the minimal version of this right approximation.

  \medskip
  For a vertex $i\in Q_0$, let
  $P(i)$ denote the $k$-linear projective representation of the quiver $Q$
  which corresponds to the vertex $i$, it gives rise to a functor
  $$P_i:\mod\Lambda\to \mod \Lambda Q, \;A\mapsto A\otimes_k P(i).$$
%  Similarly, the representation $\Rad P(i)$ give rise to the functor
%  $$P'_i:\mod \Lambda\to \mod\Lambda Q, A\mapsto A\otimes_k\Rad P(i).$$
%  Thus, $P'_\ell(A)$ is the $\Lambda$-representation of $Q$ which, for each
%  non-trivial path $\ell\to i$
%  has a copy of $A$ in position $i$.

  \medskip
  Suppose $M\in\mod\Lambda Q$ is a module and $\alpha:v\to w$ an arrow.
  We construct a module $N=\Mono_\alpha(M)$ which is such that
  the map $N_\alpha$ is monic and the sum
  $\Ima N_\alpha\oplus\sum_{t(\beta)=w,\beta\neq\alpha}\Ima N_\beta$ is direct.

  \medskip
  Choose an injective envelope $e:M_v\to J=E(M_v)$.
  Then $\Mono_\alpha(M)$ is given as the middle term of the short exact sequence
  $$\textstyle 0\to  P_w(J) \to \Mono_\alpha(M) \to M \to 0$$
  which is split exact in each component.
  For $\beta:i\to j$ an arrow, the map $\Mono_\alpha(M)_\beta$ is given by the vertical map
  in the middle of the diagram
  $$
  \begin{tikzcd}[ampersand replacement=\&]
    0 \arrow[r]
    \& \bigoplus\limits_{p:w\to i}J
    \arrow[r,"{\left(\begin{smallmatrix}0 \\ 1\end{smallmatrix}\right)}"]
    \arrow[d, "{\rm incl}"]
    \& M_i\oplus\bigoplus\limits_{p:w\to i}J
    \arrow[r, "{\left(\begin{smallmatrix} 1 & 0\end{smallmatrix}\right)}"]
    \arrow[d, "{\left(\begin{smallmatrix} M_\beta & 0\\E_\beta&{\rm incl}\end{smallmatrix}\right)}"]
    \& M_i \arrow[r]  \arrow[d, "M_\beta"] \& 0\\
    0 \arrow[r]
    \& \bigoplus\limits_{p:w\to j}J
    \arrow[r, "{\left(\begin{smallmatrix}0 \\ 1\end{smallmatrix}\right)}"]
    \& M_j\oplus \bigoplus\limits_{p:w\to j}J
    \arrow[r, "{\left(\begin{smallmatrix}1 & 0\end{smallmatrix}\right)}"]
    \& M_j \arrow [r] \& 0 \\
  \end{tikzcd}
  $$
  where $E_\beta=0$  unless $\beta=\alpha$ in which case
  $$\textstyle E_\alpha={\;e\; \choose 0}:M_v\to J\oplus \bigoplus\limits_{\substack{p:v\to j\\
      p\neq\alpha}} J,$$
  where the first component of the module on the right
  is the summand corresponding to the path $p=\alpha$.

  \medskip
  The following two results are straightforward.

  \begin{lem}
    \label{lemma-sepmon}
    Let $M\in\mod \Lambda Q$ be a module and $\alpha,\beta$ arrows.
    \begin{enumerate}
    \item The $\Lambda Q$-module $N=\Mono_\alpha(M)$ satisfies that
      $N_\alpha$ is monic and the sum
      $\Ima N_\alpha\oplus\sum_{t(\beta)=w,\beta\neq\alpha}\Ima N_\beta$ is a direct sum.
    \item
      Suppose the arrows $\alpha$, $\beta$ are such that there is no path
      from $t(\alpha)$ to $s(\beta)$, or from $t(\beta)$ to $s(\alpha)$. Then

      \smallskip
      \centerline{$\hfill\hfill \Mono_\alpha(\Mono_\beta(M))\cong \Mono_\beta(\Mono_\alpha(M)).$
        \hfill        $\s$}
      \end{enumerate}
  \end{lem}

  Note that if there are paths
      from $t(\alpha)$ to $s(\beta)$, then
     for each path, $\Mono_\beta(\Mono_\alpha(M))$ has an additional copy of $E(M_{s(\alpha)})$
  in the $t(\beta)$-branch  and in each branch corresponding to a successor of $t(\beta)$. 

  \medskip
  For a list of arrows, $A=(\alpha_s,\ldots,\alpha_1)$ such that
  for $i>j$ there is no path from $t(\alpha_i)$ to $s(\alpha_j)$,
  we define
  $$\Mono_A(M) \;=\; \Mono_{\alpha_1}(\ldots\Mono_{\alpha_s}(M)\ldots).$$
  In particular, if $v$ is a vertex in $Q$ and $(\alpha_s,\ldots,\alpha_1)$
  the list of all arrows ending in $v$, we put
  $$\Mono_v(M) \;=\; \Mono_{\alpha_1}(\ldots\Mono_{\alpha_s}(M)\ldots).$$
  Since $Q$ is an acyclic quiver, we can always arrange the set of all arrows such that
  for $i>j$ there is no path from $t(\alpha_i)$ to $s(\alpha_j)$.
   For such an arrangement of arrows, we put
  $$\Mono(M) \;=\; \Mono_{\alpha_1}(\ldots\Mono_{\alpha_s}(M)\ldots),$$
  where $s$ is the number of arrows in $Q$.

      \begin{prop}
      \label{prop-rapprox}
      Let $M\in\mod\Lambda Q$ be a module.
      \begin{enumerate}
      \item Let $A=(\alpha_s,\ldots,\alpha_1)$ be a sequence of arrows such that
        for $i>j$ there is no path from $t(\alpha_i)$ to $s(\alpha_j)$.
        If $\Mono_A(M)$ is separated monic, then the epimorphism
        in the short exact sequence
        $$\qquad 0\longrightarrow \bigoplus_{i=1}^s P_{t(\alpha_i)}(E(M_{s(\alpha_i)}))
        \longrightarrow\Mono_A(M)\longrightarrow M\longrightarrow 0 $$
        is a right approximation for $M$ in $\Gproj\Lambda Q$.
      \item In particular, if $A$ is the sequence of arrows
        in the definition of $\Mono(M)$, the epimorphism
        in the short exact sequence
        $$\qquad0\longrightarrow \bigoplus_{i=1}^s P_{t(\alpha_i)}(E(M_{s(\alpha_i)}))
        \longrightarrow\Mono(M)\longrightarrow M\longrightarrow 0 $$
        is a right approximation for $M$ in $\Gproj\Lambda Q$.
      \end{enumerate}
    \end{prop}
    \begin{proof}
      See \cite{EJ2},Theorem 11.5.1.
    \end{proof}

  \medskip
  \begin{defin} For $M\in\mod\Lambda Q$, the {\it minimal monomorphism}
    $\Mimo(M)$ is given by the right minimal version $\Mimo(M)\to M$
    of the map $\Mono(M)\to M$.
  \end{defin}

  \medskip
  Since  $\Mono(M)$ is a finite length module (Lemma~\ref{prop-rapprox}),
  the approximation $\Mono(M)\to M$ has a right minimal version;
  this map $\Mimo(M)\to M$ is
  the minimal right approximation for $M$ in $\Gproj\Lambda Q$.
  To obtain the right minimal version of $\Mono(M)\to M$,
  we need to split off some projective direct summand
  from $\Mono(M)$.

  \medskip
  We have

  \begin{prop} \label{split}
    Let $Q$ be a finite acyclic quiver, $n$ a source in $Q$, $\Lambda$ a selfinjective
    algebra and $M$ a separated monic representation of $Q$ over $\Lambda$.
    If  $M_n=I_n\oplus M'_n$ with $I_n$ being a projective $\Lambda$-module,
    then $P_n(I_n)$ is a direct summand of $M$.
\end{prop}

\begin{proof}
  There is an embedding $\sigma: I_n\to M_n$ which induces
a monomorphism $P_n(\sigma): P_n(I_n)\to M$ since $M$ is separated monic
(Lemma 2.3 in \cite{LZ1}). With $S=\Cok P_n(\sigma)$,
we have an exact sequence
$$0\to P_n(I_n)\stackrel{ P_n(\sigma)}\to M \stackrel{ \pi}\to S \to 0.$$
Since $P_n(I_n)$ is a projective-injective object in $\Gproj \Lambda Q $, then to prove the assertion, we just need to prove that $S$ is separated monic.

\medskip

For each $i\in Q_0\setminus \{n\}$, the above short exact sequence
gives the commutative diagram:

\[ \begin{tikzcd}[column sep=small]
  0\arrow{r}{}& P_n(I_n)_{i^-}\arrow{r}{   }\arrow{d}{u_{P,i}} & M_{i^-}
  \arrow{d}{u_{M,i}} \arrow {r}{ } &  S_{i^-}\arrow{d}{u_{S,i}} \arrow{r}{ }
                                                           &  0 \\
  0 \arrow{r}{} &P_n(I_n)_i\arrow{r}{ } & M_i\arrow{r}{ } & S_i  \arrow{r}{} & 0 \\
\end{tikzcd}
\]
Recall that a module of the form $P_w(V)$ where $w\in Q_0$ and $V\in\mod\Lambda$
has the property that $u_{P_w(V),i}$ is an isomorphism for every vertex $i\neq w$.
Here, since $i\neq n$, the left hand map always is an isomorphism.

Hence, by the Snake Lemma, for each $i\in Q_0$, $$u_{S,i}: S_{i^-}\to S_i$$
is a monomorphism which follows from the facts that $M$ is separated monic
and $\Cok (u_{P,i})=0$.
\end{proof}

\begin{prop} \label{nonprojective} Let $Q$ be a finite acyclic quiver,
  $\Lambda$ a selfinjective algebra
    and $M$ be a separated monic representation of $Q$ over $\Lambda$.
    There is an exact sequence
    $$0\to I\to M\to N\to 0$$ such that $I$ is a projective object in
    $\Gproj \Lambda Q$ and no branch of the $\Lambda Q$-module
    $N$ has a non-zero projective direct summand.
\end{prop}
\begin{proof}
  Use induction on the number $n=|Q_0|$ of vertices in $Q$.
  If  $n=1$, then the assertion is clear.
  Suppose that the assertion is true for quivers with fewer vertices.

Let $n$ be a source of $Q$.
Suppose that $M_n=I_n\oplus L_n$ such that $I_n$ is a projective $\Lambda$-module and $L_n$ has no non-zero projective direct summand. By Proposition \ref{split}, we have an exact sequence

$$0\to P_n(I_n)\stackrel{ P_n(\sigma)}\to M \stackrel{ \pi}\to S \to 0,$$
with $S= \Cok   P_n(\sigma)$ being a separated monic representation of $Q$
over $\Lambda$ and $S_n=L_n$.
Restricting $S$ to the full subquiver $Q'$ with set of vertices $Q_0\setminus \{n\}$,
we get a separated monic representation $ S'$ of $Q'$ over $\Lambda$.
Since $|Q'_0|=n-1$,  by the hypothesis, we have an exact sequence
$$0\to I' \to  S'  \to N'\to 0,$$
such that  $I'$ is a projective object in $\Gproj \Lambda Q'$ and
no branch of $N'$ has a non-zero projective direct summand.
Consider $I'$ as a representation of $Q$;
let $N$ be the target of the cokernel map in the diagram
\[ \begin{tikzcd}[column sep=small]
  0\arrow{r}{}& I'\arrow{r}{   }\arrow[d,equal] &S'
  \arrow[d,hook] \arrow {r}{ } &  N'\arrow[d,dotted] \arrow{r}{ }
                                                           &  0 \\
  0 \arrow{r}{} &I'\arrow{r}{ } & S\arrow{r}{ } & N  \arrow{r}{} & 0 \\
\end{tikzcd}
\]
 
  Thus, $N_i=N_i'$ for $i\neq n$ and $N_n=S_n$.
  In particular, no branch of $N$ contains a non-zero projective-injective direct summand. Take $I$ to be the pullback of $ \pi$ and $I'\to S$, given by the following
commutative diagram

\[ \begin{tikzcd}
  &&0\arrow{d}&0\arrow{d}&\\
0 \arrow{r}{} & P_n(I_n) \arrow{r}{ } \arrow[d,equal] & I\arrow{d}{ } \arrow{r}{ }
                                                           &I'\arrow{r}{} \arrow[d] & 0 \\
  0 \arrow{r}{} & P_n(I_n)\arrow{r}{ } &M\arrow[r,"\pi"] \arrow{d}{ }& S\arrow{r}{}\arrow{d}{ } & 0 \\
&&N\arrow{d}\arrow[r,equal]&N\arrow{d}&\\
 &&0&0 &\\
\end{tikzcd}
\]
with exact rows and columns.
Since $ P_n(I_n) $ and $I'$ are projective objects in $\Gproj \Lambda Q$,
so is $I$.
\end{proof}

\medskip
If $M$ is a separated monic representation of $Q$ over $\Lambda$ and
for each vertex $i\in Q_0$, $M_i$ is a projective $\Lambda$-module, then by
Proposition \ref{nonprojective}, we obtain an exact sequence $$0\to I\to M\to N\to 0$$
such that $I$ is a projective object in $\Gproj \Lambda Q$ and no branch
of $N$ has a non-zero projective direct summand.  It follows that  $N=0$ and
$M\cong I$. We have shown the following result:

\begin{cor} \label{pro} Let $Q$ be a finite acyclic quiver,
  $\Lambda$ a selfinjective algebra
    and $M$ be a separated monic representation of $Q$ over $\Lambda$.
    If $M_i$ is a projective $\Lambda$-module for each vertex $i\in Q_0$,
    then $M$ is a projective representation $Q$ over $\Lambda$.\qed
\end{cor}

  %--------------------------------------------------------------------------------
  \subsection{Lifting isomorphisms from the stable category}
  %--------------------------------------------------------------------------------

  The following result adapts \cite[Theorem 4.2]{RS} to our situation.

  \begin{prop}\label{Mono(XY)}
    Let $X,Y$ be representations of $Q$  over  $\Lambda$ such that none of the components $X_i$, $Y_i$
    has a nonzero injective direct summand.
    Then $\Mono(X)\cong\Mono(Y)$ if and only if for each vertex $i$, there
    is an isomorphism $f_i:X_i\to Y_i$ such that for each arrow $\alpha:i\to j$, the map
    $Y_\alpha f_i-f_j X_\alpha$ factors through an injective object.

    \begin{center}\begin{tikzcd}[row sep=small, column sep=small]
      X_i \arrow[r,"f_i"] \arrow[d,swap,"X_\alpha"] & Y_i \arrow[d,"Y_\alpha"] \\ X_j \arrow[r,swap, "f_j"] & Y_j
    \end{tikzcd}\end{center}

  \end{prop}

  \begin{proof}
    Assume first that $g:\Mono(X)\to \Mono(Y)$ is an isomorphism. Write $\Mono(X)_i=X_i\oplus D_i$,
    $\Mono(Y)_i=Y_i\oplus E_i$ with $D_i$ and $E_i$ injective $\Lambda$-modules,
    and let $g_i=\begin{smallpmatrix} a_i & b_i\\ c_i & d_i\end{smallpmatrix}$.
    Since $X_i$ and $E_i$, and $Y_i$ and $D_i$ have no isomorphic indecomposable direct summands,
    the maps $a_i$ and $d_i$ must be isomorphisms.

    For each arrow $\alpha:i\to j$, we have $\Mono(Y)_\alpha\circ g_i=g_j \circ\Mono(X)_\alpha$,
    hence, by restriction, there is the commutative diagram
    \begin{center}\begin{tikzcd}[ampersand replacement=\&,  column sep=small]
        X_i \arrow[r, "{\begin{smallpmatrix} a_i \\ c_i \end{smallpmatrix}}"]
          \arrow[d, swap, "{\begin{smallpmatrix} X_\alpha\\ u\end{smallpmatrix}}"]
        \& \Mono(Y)_i \arrow[d,"{\begin{smallpmatrix}Y_\alpha & v\end{smallpmatrix}}"] \\
        \Mono(X)_j \arrow[r,swap,"{\begin{smallpmatrix}a_j & b_j\end{smallpmatrix}}"] \& Y_j
    \end{tikzcd}\end{center}
    and $Y_\alpha a_i-a_j X_\alpha=b_j u-vc_i$ factors through an injective module.

    \smallskip
    For the converse, assume isomorphisms $f_i:X_i\to Y_i$ are given such that
    for each arrow $\alpha:i\to j$, the map
    $Y_\alpha f_i-f_j X_\alpha$ factors through an injective object.
      We recall that $\Mono(M)=\Mono_{\alpha_1}(\ldots\Mono_{\alpha_s}(M)\ldots)$
      and proceed along the sequence of arrows $\alpha_s,\alpha_{s-1},\ldots,\alpha_1$.

      \smallskip
      Let $k\in\{s,\ldots,1\}$ and write
      $U=\Mono_{\alpha_{k+1}}(\ldots\Mono_{\alpha_s}(X)\ldots)$ and
      $V=\Mono_{\alpha_{k+1}}(\ldots\Mono_{\alpha_s}(Y)\ldots)$.
      Suppose we have already constructed a map (not necessarily a homomorphism)
      $g^{k+1}:U\to V$ which is a $\Lambda$-isomorphism in each branch and which is 
      such that each square corresponding to $\alpha_{k+1},\ldots,\alpha_s$
      is commutative.  Write $\alpha=\alpha_k$, say $\alpha:i\to j$, and consider
      the square corresponding to this arrow, recall that it may not be commutative.
      $$\begin{tikzcd}[ampersand replacement=\&]
        X_i \arrow{r}{f_i} \arrow[d,swap,"{X_\alpha\choose0}"] \& Y_i \arrow{d}{{Y_\alpha\choose0}} \\
        X_j\oplus Q_j
          \arrow[r,swap,"{\begin{smallpmatrix}f_j & a_j\\ 0\; & b_j\end{smallpmatrix}}"] \&
          Y_j\oplus R_j\\
      \end{tikzcd}
      $$
      (The sequence of arrows is chosen such that we have $U_i=X_i$, $V_i=Y_i$
      and $g^{k+1}_i=f_i$, but in the $j$-component there may be additional
      projective summands $Q_j$, $R_j$ with $b_j:Q_j\to R_j$ an isomorphism.)
      In the construction $\Mono_\alpha$ for $U$ and $V$, we have added injective
      envelopes $d_i:X_i\to I_i$, $e_i:Y_i\to J_i$.
      The given isomorphism $f_i:X_i\to Y_i$ can be extended to an isomorphism
      $u_i:I_i\to J_i$ such that $u_i d_i=e_i f_i$.
    The map $f_jX_\alpha-Y_\alpha f_i$ factors through an injective module,
    hence factors through the injective envelope $d_i:X_i\to  I_i    $, so there is
    $h_\alpha: I_i  \to Y_j$ such that $f_jX_\alpha-Y_\alpha f_i=-h_\alpha d_i$.
    We construct
      the map $g^k:\Mono_\alpha(U)\to \Mono_\alpha(V)$. The square corresponding to $\alpha$
      is as follows.
      $$\begin{tikzcd}[ampersand replacement=\&]
        X_i \arrow{r}{f_i}
        \arrow[d,swap,"{\begin{smallpmatrix}X_\alpha\\0\\d_i\end{smallpmatrix}}"]
          \& Y_i \arrow{d}{\begin{smallpmatrix}Y_\alpha\\0\\e_i\end{smallpmatrix}} \\
        X_j\oplus Q_j\oplus I_i
        \arrow[r,swap,"{\begin{smallpmatrix}f_j & a_j & h_\alpha \\ 0\; & b_j & \;0\\
              0\; & 0 &\; u_i \end{smallpmatrix}}"] \&
          Y_j\oplus R_j\oplus J_i\\
      \end{tikzcd}
      $$
       
  For each remaining vertex $\ell\neq i,j$ we construct the map
  $g^k_\ell:\Mono_\alpha(U)_\ell\to \Mono_\alpha(V)_\ell$ as follows.
  Note that $\Mono_\alpha(U)_\ell=U_\ell\oplus P_j(I_i)_\ell$ where $P_j(I_i)_\ell=\bigoplus_{p:j\to \ell} I_i^{(p)}$,
  the sum being taken over all paths $p$ from $j$ to $\ell$ and each summand $I_i^{(p)}=I_i$ being marked
  by the path to which it corresponds.
  Define
  $$g_\ell^k=\left(\begin{array}{cc}g_\ell^{k+1} & (V_p\circ {h_\alpha\choose 0})_p\\ 0 & \bigoplus u_i^{(p)}\end{array}\right):
    U_\ell\oplus\bigoplus I_i^{(p)} \longrightarrow V_\ell\oplus \bigoplus J_i^{(p)};$$
    here, ${h_\alpha\choose 0}:I_i\to V_j=Y_j\oplus R_j$ is part of the map $g^k_j$ pictured
    above, $V_p=V_{\beta_t}\circ\cdots\circ V_{\beta_1}$ is the composition of the
    structural maps in $V$ if the path is given as $p=\beta_t\cdots\beta_1$,
    and $\bigoplus_{p:j\to\ell} u_i^{(p)}:\bigoplus I_i^{(p)}\to\bigoplus J_i^{(p)}$ is the diagonal map
    with entries $u_i^{(p)}=u_i$.

    \smallskip
    The map $g^k:\Mono_\alpha(U)\to \Mono_\alpha(V)$ thus constructed has the following properties:
    Like $g^{k+1}$, for each arrow the corresponding square is commutative, up to a morphism which
    factors through an injective object.  Like $g^{k+1}$, for each arrow $\alpha_{k+1},\ldots,\alpha_s$,
    the corresponding square is commutative.  And in addition, the square for the arrow $\alpha=\alpha_k$
    is commutative.  Hence we can continue our iteration.
     
        \smallskip
        With the iteration completed, it follows that the map $g=g^1:\Mono(X)\to\Mono(Y)$
        is a homomorphism since for each
    arrow, the corresponding diagram is commutative, and hence an isomorphism
    since each branch map is.

    \smallskip
    We illustrate this map in the
    example of the quiver $Q:{\sssize 1}\stackrel\alpha\to {\sssize2} \stackrel\beta\to {\sssize3}$. For representations $X,Y$ of $Q$ over $\Lambda$ and isomorphisms $f_i:X_i\to Y_i$
      as in Proposition~\ref{Mono(XY)},
      the corresponding map $g:\Mono_\alpha(\Mono_\beta(X))\to \Mono_\alpha(\Mono_\beta(Y))$
      is as follows.

    \begin{center}\begin{tikzcd}[ampersand replacement=\&,  column sep=large, row sep=large]
        X_1 \arrow[r,"f_1"] \arrow[d,swap,"\begin{smallpmatrix}X_\alpha\\ d_1\end{smallpmatrix}"]
        \& Y_1 \arrow[d,"\begin{smallpmatrix}Y_\alpha\\ e_1\end{smallpmatrix}"] \\
        X_2\oplus  I_1
      \arrow[r,"{\begin{smallpmatrix}f_2 & h_\alpha\\0 & u_1\end{smallpmatrix}}"]
        \arrow[d,swap,"{\begin{smallpmatrix}X_\beta & 0\\ d_2 & 0 \\ 0 & 1\end{smallpmatrix}}"]
        \& Y_2\oplus  J_1   \arrow[d,"{\begin{smallpmatrix}Y_\beta& 0\\e_2 & 0\\ 0 & 1\end{smallpmatrix}}"] \\
        X_3\oplus  I_2
     \oplus  I_1
        \arrow[r,swap,"{\begin{smallpmatrix}f_3 & h_\beta & Y_\beta h_\alpha\\ 0 & u_2 & e_2h_\alpha\\ 0 & 0 & u_1\end{smallpmatrix}}"]
        \& Y_3\oplus  J_2
    \oplus  J_1  \end{tikzcd}\end{center}
  \end{proof}

  \begin{rem} From the above proof we can see that the condition
    that none of the components $X_i$, $Y_i$
    have a nonzero injective direct summand is not necessary for the sufficiency part
    of the proposition.
  \end{rem}

%================================================================================
\section{The equivalence of stable Gorenstein-projective categories}
%================================================================================
\label{section-equivalence}

Let $Q$ be a finite acyclic quiver with sink $v$ and
$\Lambda$ a selfinjective algebra.
The quiver $Q(v)$ is obtained from $Q$ by reversing all the arrows
ending in $v$.
In this section, we
  construct the equivalence between $\uGproj \Lambda Q$ and
  $\uGproj \Lambda Q(v)$ from Theorem~\ref{theorem-intro1}
  on the level of objects and morphisms.
%--------------------------------------------------------------------------------
\subsection{Reflection on objects}
%--------------------------------------------------------------------------------
\label{sec-reflection-objects}

\medskip

Let $X\in \Gproj \Lambda Q$.
The $\Lambda$-homomorphism $u_X:X_{v^-}\to X_v$ from
Section~\ref{section-monic-notation} and a projective cover
$\delta_X:C_X\to X_{\bar v}$ of the cokernel of $u_X$
give rise to the following commutative
diagram with exact rows.
\[ \begin{tikzcd}
  0 \arrow{r}{} & \Omega_X \arrow{r}{\sigma_X} \arrow[swap]{d}{k_X} & C_X\arrow{d}{ t_X} \arrow{r}{\delta_X}
                                                           & X_{\bar v}\arrow{r}{} \arrow[d, equal]& 0 \\
  0 \arrow{r}{} & X_{v^-}\arrow{r}{u_X} & X_v\arrow{r}{\pi_X} & X_{\bar v} \arrow{r}{} & 0 \\
\end{tikzcd}
\]

We index the components of the kernel map $k_X$ using the arrows
$\beta$ in $Q(v)$ starting at $v$, thus
$k_X=(k_{X,\beta})_{s(\beta)=v}$ for $\Lambda$-homomorphisms
$k_{X,\beta}:\Omega_X\to X_{t(\beta)}$.

\medskip
We define the $\Lambda Q(v)$-module $X'(v)$ in terms of $\Lambda$-modules
$$X'(v)_i=\left\{\begin{array}{ll} X_i & \text{if}\;i\neq v\\ \Omega_X & \text{if}\; i=v\end{array}\right.$$
and $\Lambda$-linear morphisms $X'(v)_\beta$ where $\beta$ is an arrow in $Q(v)$
$$X'(v)_\beta=\left\{\begin{array}{ll} X_\beta & \text{if}\;
 s(\beta)\neq v\\
k_{X,\beta} & \text{if}\;s(\beta)=v.\end{array}\right.$$

\medskip

 Note that the representation $X'(v)$ may not be
  separated monic.

 Hence we define $$X(v)=\Mono_v(X'(v)).$$
    It follows that $X(v)$ is separated monic, hence
    Gorenstein-projective, and the epimorphism in the
    short exact sequence
    $$0\to \bigoplus_{s(\beta)=v} P_{t(\beta)}(C_X)\to X(v)\to X'(v)\to 0$$
    in  Proposition~\ref{prop-rapprox}~(1)
    is a right approximation for $X'(v)$ in $\Gproj\Lambda Q$.
    In particular, we have:

    \begin{lem}\label{unique}
      The module $X(v)$ is determined uniquely, up to isomorphy.
    \end{lem}

For later use we summarize the construction for $X(v)$:

    \smallskip
    Let $X\in \Gproj\Lambda Q$ and $v\in Q_0$ a sink.
    Choose as above a projective cover $C_X\to X_{\bar v}$ with kernel
    $\sigma_X:\Omega_X\to C_X$ and maps
    $k_{X,\beta}:\Omega_X\to X_{t(\beta)}$, where $\beta:v\to t(\beta)$
    is an arrow in $Q(v)$ starting at $v$.
  Let
$$P=\bigoplus_{s(\beta)=v} P_{t(\beta)}(C_X)$$
  be the direct sum, indexed by the arrows starting in $v$,
  of the projective $\Lambda Q(v)$-modules
  $P_{t(\beta)}(C_X)=C_X\otimes_k P(t(\beta))$
 where $P(t(\beta))$ is the
   projective $kQ(v)$-module
 corresponding to the vertex $t(\beta)$.
  Then $X(v)$ is as follows:
  $$X(v)_i=\left\{\begin{array}{ll} (X\oplus P)_i & \text{if}\;i\neq v\\
  \Omega_X & \text{if}\; i=v\end{array}\right.$$
  and
  $$X(v)_\beta=\left\{\begin{array}{ll} (X\oplus P)_\beta
  & \text{if}\;   s(\beta)\neq v\\
  \left( \begin{smallmatrix}k_{X,\beta}\\ \sigma_X\\ 0 \end{smallmatrix}\right) &
  \text{if}\;s(\beta)=v.\end{array}\right.$$
  where the components of the target of the map 
  $$
  \left( \begin{smallmatrix}k_{X,\beta}\\ \sigma_X\\ 0 \end{smallmatrix}\right):
  \quad\Omega_X\longrightarrow X_{t(\beta)}\oplus C_{X,\beta}\oplus
  \bigoplus\limits_{ \begin{smallmatrix}p:v\to t(\beta)\\p\neq \beta \end{smallmatrix}}
  C_{X, p}
  $$
  are indexed by the paths $p:v\to t(\beta)$ in the construction
  of $P$.  Here, $C_{X,\beta}=C_X$ for the given arrow $\beta$ starting at $v$
  and  $C_{X,p}=C_X$ for each path $p\neq\beta$ from $v$ to $t(\beta)$,
  but the only nonzero map corresponds to the component of $P$ given by
  the arrow $\beta$.

%--------------------------------------------------------------------------------
\subsection{The converse construction and density}
%--------------------------------------------------------------------------------
\label{sec-converse}

Now, consider the converse of the above construction defined in Subsection 3.1. Let $Y$ be an object in $\Gproj \Lambda Q(v)$. Then the injective envelope $\sigma_Y:Y_v\to C$ and the homomorphism
$k: Y_v\to Y_{v^+}=\bigoplus\limits_{s(\beta)=v} Y_{t(\beta)}$ give a pushout diagram

\[ \begin{tikzcd}
  0 \arrow{r}{} &Y_v \arrow{r}{\sigma_Y} \arrow[d, "k" ] & C \arrow[d] \arrow{r}{ }
   & \Cok {\sigma_Y} \arrow{r}{} \arrow[d, equal] & 0 \\
  0 \arrow{r}{} & Y_{v^+}\arrow[r, "b"] & B\arrow{r}{ } &\Cok {\sigma_Y} \arrow{r}{} & 0 \\
\end{tikzcd}
\]
Write $b=(b_\alpha)$ where $\alpha^{op}$ runs over all arrows in $Q(v)$ starting at $v$.
Let $X=(X_i, X_\alpha)$ be the representation of $Q$ over $\Lambda$, where
 \begin{center}
$X_i=\left\{\begin{array}{ll} B, & \text{if}\; i=v; \\
Y_i ,& \text{otherwise}. \end{array}\right.$ and
 $X_\alpha=\left\{\begin{array}{ll} b_\alpha, & \text{if}\; t(\alpha)=v; \\
Y_\alpha, & \text{otherwise}. \end{array}\right.$
\end{center}
 Since $b$ is a monomorphism and $Y$ is separated monic, then $X$
 is a separated monic representation of $Q$ over $\Lambda$.
 Moreover, $X_{\bar v}=\Cok \sigma_Y $.
 Note that $C\to \Cok\sigma_Y$ is a projective cover. It follows that
 $X'(v)\cong Y$ in $\Gproj \Lambda Q(v)$ (see Subsection~\ref{unique}).
 Since $X'(v)$ is in $\Gproj\Lambda Q(v)$ and projective objects are injective in  $\Gproj\Lambda Q(v)$
 (see Remark 10.2.2 in \cite{EJ2} ),
 the exact short sequence in Subsection 3.1 splits.
 Hence $X(v)\cong X'(v)\oplus P$, so $X(v)\cong Y$ in $\underline{\Gproj}\Lambda Q(v)$ and we have
 shown that the construction is dense.

%--------------------------------------------------------------------------------
\subsection{Reflection on homomorphisms}
%--------------------------------------------------------------------------------
\label{sec-refl-homom}

Let $f:X\to Y$ be a homomorphism in $\Gproj \Lambda Q $.
We define a corresponding homomorphism
$f(v):X(v)\to Y(v)$ in $\Gproj\Lambda Q(v)$.
Let $u_X$, $u_Y$ be the $\Lambda$-homomorphisms used in the construction
for $X(v)$ and $Y(v)$, they give rise to the commutative diagram
$$\begin{tikzcd}[row sep=small, column sep=small]
  X_{v^-}\arrow[r,"u_X"]\arrow[d,swap,"f_{v^-}"] & X_v\arrow[d,"f_v"]\\
  Y_{v^-}\arrow[r,swap,"u_Y"] & Y_v \end{tikzcd}
$$
where $f_{v^-}=\bigoplus\limits_{w\in v^-} f_w$ is the componentwise map.  Hence
there is a homomorphism  between the cokernels,
  $f_{\bar v}: X_{\bar v}\to Y_{\bar v}$
such that  $f_{\bar v}\pi_X=\pi_Yf_{v}$
as in the diagram below.
Fix projective covers $\delta_X: C_X\to X_{\bar v}$ and
$\delta_Y: C_Y\to Y_{\bar v}$ of $X_{\bar v}$ and $Y_{\bar v}$, respectively.
Let $\Omega_X=\ker \delta_X$, $\Omega_Y=\ker \delta_Y $,
then there are homomorphisms $t_X:C_X\to X_v$ and $k_X: \Omega_X\to X_{v^-}$
such that $\pi_Xt_X=\delta_X, u_Xk_X=t_X\sigma_X$ and $t_Y:C_Y\to Y_v$ and
$k_Y: \Omega_Y\to Y_{v^-}$ such that $\pi_Yt_Y=\delta_Y, u_Yk_Y=t_Y\sigma_Y$.
Since $\delta_Y$ is an epimorphism, there exists $C_f: C_X\to C_Y$
such that $\delta_YC_f=f_{\bar v} \delta_X$ and $\Omega_f:\Omega_X\to\Omega_Y $
such that $\sigma_Y\Omega_f=C_f \sigma_X$.
So we have a diagram in which every square is commutative except
for the left and the middle ones.
\[\begin{tikzcd}[row sep=small, column sep=small]
  0\arrow[rr]& &
  \Omega_X\arrow[dr,near end, "\Omega_f"]\arrow{dd}[near start]{k_X}
  \arrow[rr, "\sigma_X"]& &
  C_X\arrow[dddl,"\omega_f ",dotted] \arrow[dr,near end, "C_f"]
  \arrow{dd}[near start]{t_X} \arrow[rr, "\delta_X"]& &
 X_{\bar v}\arrow[dr, "f_{\bar v}"]\arrow[dd,  equal] \arrow[rr]& &0\\
 &0\arrow[rr, crossing over]& & \Omega_Y
 \arrow[rr,   crossing over, "\sigma_Y\ \ \ \ \ \ \ "]& &
 C_Y \arrow[rr, crossing over, "\delta_Y\ \  \ \ \ \ \ "]& &
 Y_{\bar v} \arrow[rr]& &0\\
 0\arrow[rr, crossing over]& &  X_{v^-}\arrow[dr,"f_{v^-}"]
 \arrow[rr, "u_X\ \  \ \ \ \ \ "]& &
 X_v\arrow[dr, "f_v"] \arrow[rr,  "\pi_X\ \  \ \ \ \ \ "]& &
 X_{\bar v}\arrow[dr, "f_{\bar v}"] \arrow[rr]& &0\\
 & 0\arrow[rr]& &
 Y_{v^-} \arrow[from=uu, "k_Y", near start, crossing over] \arrow[rr, "u_Y"]& &
 Y_v \arrow[from=uu, "t_Y", near start, crossing over] \arrow[rr, "\pi_Y"]& &
 Y_{\bar v} \arrow[from=uu, equal, crossing over] \arrow[rr]& &0\\
\end{tikzcd}.\]
Since $\pi_Y(t_YC_f-f_v t_X)=0$, then there is a morphism
$\omega_f: C_X\to Y_{v^-}$ such that $u_Y \omega_f=t_YC_f-f_v t_X$.
And $\omega_f\sigma_X=k_Y\Omega_f-f_{v^-}k_X$ follows from the
injectivity of $u_Y$. We index $\omega_f=(\omega_{f,\beta})_{\{\beta: s(\beta)=v\}}.$
This is to say, for each $\beta: v\to i $ in $ Q(v)_1$, we have
$\omega_{f,\beta}\sigma_X=k_{Y,\beta}\Omega_f-f_ik_{X,\beta}$.

 \begin{rem}
   Note that the map $f'(v):X'(v)\to Y'(v)$ defined by
   $$f'(v)_i=\left\{\begin{smallmatrix}\Omega_f, & \text{if}\;i=v\\
   f_i,&\text{if}\;i\neq v\end{smallmatrix}\right.$$
   is not necessarily a homomorphism since the left square in the above diagram may not commute.
 \end{rem}

Let $f(v)=(f(v)_i)$ where $f(v)_v= \Omega_f $ and for each $i\neq v$
$$f(v)_i={\left(\begin{smallmatrix} f_i &  (Y_{p'} \omega_{f,\beta})_{p} \\
    0 &C_fE  \end{smallmatrix} \right)}:
X_i\oplus\bigoplus\limits_{p=p'\beta} C_{X,p}\to Y_i\oplus
\bigoplus\limits_{p=p'\beta} C_{Y,p}, $$
where $p$ runs over all non-trival paths from $v$ to $i$ in $Q(v)$
 such that $\beta$ is the starting arrow of $p$,
  $p'$ is the subpath of $p$ such that $p=p'\beta$
  and $E$ is the unit matrix of the size $|\{p: v\to i\}|$.
  Note that when $p$ is an arrow, then $p'$ is a trival path $t(p)$.
  In this case $Y_{p'}$ is an identity morphism of $Y_{t(p)}$.
  One can check that $f(v): X(v) \to Y(v)$ is a homomorphism in
  $\Gproj \Lambda Q(v)$.

  \medskip
  We illustrate this construction in an example.

  \begin{ex}
    Let $Q$ be the quiver $1\to 2\to 3$, then $Q(3)$ is the quiver
    $1\to 2\leftarrow 3$. Let $k$ be a field, $A=k[x]/(x^2)$  the
    algebra of dual numbers, $\sigma: k\to A$ an embedding and
    $\pi: A\to k$ the canonical epimorphism.
    Consider the homomorphism $f:X\to Y$ between two separated monic
    representations given by the diagram.
    $$\begin{tikzcd}[row sep=small]
      X: & 0 \arrow[r] \arrow[d] & k \arrow[r,"\sigma"] \arrow[d,"\sigma"]
      & A \arrow[d,"1_A"]\\
      Y: & k \arrow[r,swap,"\sigma"] & A \arrow[r,swap,"1_A"] & A
    \end{tikzcd}
    $$
    Since $X_{\bar 3}=\Cok\sigma\cong k$, we have $C_X=A, \Omega_X=k$.
    Then
    $$X'(3)=0\to k \stackrel{id_k}\longleftarrow k,\qquad
    P_2(C_X)=0\to A\leftarrow 0,
    $$
    so we have
    $X(3)=0\to k\oplus A \stackrel{\left(
      \begin{smallmatrix} id_k\\ \sigma \end{smallmatrix}\right)}
    \longleftarrow k$.
    Since $Y_{\bar 3}=\Cok id_A=0$, we have $C_Y=\Omega_Y=0$.
    Then $Y(3)=k\to A \leftarrow 0$.

    \medskip
    Now, we need to build $f(3)$. We illustrate the related
    homomorphisms with the following diagram

\[\begin{tikzcd}[column sep=small, row sep=small]
  0\arrow[rr]& &
  k\arrow[dr,near end, ""]\arrow{dd}[near start]{id_k}\arrow[rr, "\sigma"]& &
  A\arrow[dddl,"\omega_f ",dotted] \arrow[dr,near end, " "]\arrow{dd}[near start]{id_X} \arrow[rr, "\pi"]& &
 k\arrow[dr, " "]\arrow[dd,  equal] \arrow[rr]& &0\\
&0\arrow[rr, crossing over]& & 0 \arrow[rr,   crossing over, " "]& &
 0 \arrow[rr, crossing over, " "]& &
 0\arrow[rr]& &0\\
 0\arrow[rr, crossing over]& & k\arrow[dr,"f_2"] \arrow[rr, "\sigma\ \  \ \ \ \ \ "]& &
 A\arrow[dr, "f_3"] \arrow[rr,  "\pi\ \  \ \ \ \ \ "]& &
 k\arrow[dr, " "] \arrow[rr]& &0\\
 & 0\arrow[rr]& &
 A \arrow[from=uu, " ", near start, crossing over] \arrow[rr, "id_A"]& &
A \arrow[from=uu, " ", near start, crossing over] \arrow[rr, " "]& &
0 \arrow[from=uu, equal, crossing over] \arrow[rr]& &0\\
\end{tikzcd}\]
where $\omega_f\sigma=-f_2 id_k=-\sigma$. Taking $\omega_f=-id_A$,
then $f(3)_2=\left(\begin{smallmatrix} \sigma, & -id_A\end{smallmatrix}\right)$.
  Then the homomorphism $f(3):X(3)\to Y(3)$ is as follows.
      $$\begin{tikzcd}[ampersand replacement=\&]
    X(3): \& 0 \arrow[r] \arrow[d]
    \& k\oplus A \arrow[d,"\left(\sigma\;\;-1_A\right)"]
      \& k \arrow[d]\arrow[l,swap,"\left({1_k}\atop\sigma\right)"] \\
      Y(3): \& k \arrow[r,swap,"\sigma"] \& A  \& 0 \arrow[l]
    \end{tikzcd}
    $$
\end{ex}

  \medskip

  In general, the construction of $f(v)$ depends on choices for
  $C_f$, $\omega_f$ and so on.
  It is straightforward to verify the following

  \begin{lem}
    The class of $f(v)$ in the stable category $\uGproj\Lambda Q(v)$
    depends only on the class of $f$ in $\uGproj\Lambda Q$.
    Moreover, the map
    $$\Hom_{\uGproj\Lambda Q}(X,Y)\to \Hom_{\uGproj\Lambda Q(v)}(X(v),Y(v)),
    \quad f\mapsto f(v)$$
    is a homomorphism of abelian groups.  \qed
  \end{lem}

%--------------------------------------------------------------------------------
\subsection{Equivalence of stable categories}
%--------------------------------------------------------------------------------

\phantom m
 The assignments $X\mapsto X(v)$,
  $f\mapsto f(v)$ studied in the previous two subsections give rise to a
  functor $F(v):\underline{ \Gproj} \Lambda Q \to \underline{ \Gproj} \Lambda Q(v)$.
  Actually, this functor $F(v)$ induces an equivalence of categories.

\medskip

\begin{thm} \label{equivaentl}
  Let $\Lambda$ be a selfinjective algebra,
  $Q$ be a finite acyclic quiver, and $v$ a sink in $Q$.
  Then the functor
  $$F(v):\underline{ \Gproj} \Lambda Q \to \underline{ \Gproj} \Lambda Q(v)$$
  yields an equivalence between the two stable categories modulo projectives.
\end{thm}

\begin{proof} In Subsections~\ref{sec-reflection-objects} and \ref{sec-refl-homom} we defined
  $F(v)$ as a construction for objects and homomorphisms.
  The proof that $F(v)$ gives rise to a functor is lengthy and but straightforward
  and is omitted. By Subsection~\ref{sec-converse}, we know that $F(v)$ is dense.

  \bigskip
Claim $1$:  $F(v)$ is faithful.
\vskip5pt
Let $f: X\to Y$ be a morphism in $ \Gproj \Lambda Q $ such that
the class of $f(v)$ is $0$ in $\underline{ \Gproj} \Lambda Q(v)$,
that is to say, $f(v)=gh$ where $X(v)\stackrel{h}\to M\stackrel{g} \to Y(v)$
and $M$ is projective in $ \Gproj  \Lambda Q(v)$.
 We decompose the branches of $X(v)$, $Y(v)$ into direct sums as in the
  summary after Lemma~\ref{unique}. 
Denote by
$$g_i = \begin{cases} g_v, &  \ i =v; \\
\left(\begin{smallmatrix}g^Y_{i} \\ (g_p)_p\end{smallmatrix}\right), &    i \neq v.\end{cases},\ h_i = \begin{cases} h_v, &  \ i =v; \\
\left(\begin{smallmatrix}h^X_{i}, & (h_p)_p\end{smallmatrix}\right), &    i \neq v.\end{cases} $$
where $g^Y_{i }: M_i\to Y_i, h^X_{i}: X_i\to M_i$, and
$g_p: M_i\to C_{Y,p}, h_p: C_{X,p}\to M_i$,
for each path $p$ from $v$ to $i$
where $p$ runs over all non-trival paths from $v$ to $i$.

\medskip
 
  We first construct a projective representation $N$ over $\Lambda Q$,
  and then define maps $h':X\to N$ and $g':N\to Y$ such that $f=g'h'$.

\medskip
Let $N=(N_i, N_{\alpha})$ be an object in $ \Gproj  \Lambda Q $ where $N_i=M_i$
 for each $i\neq v$, $N_v=(\bigoplus\limits_{t(\alpha)= v}M_{s(\alpha)})\oplus C_Y$;
 $N_\alpha=M_\alpha$ for each $\alpha$ with $t(\alpha)\neq v$, $N_\alpha: M_{s(\alpha )}\to N_v $
 are natural embedding for each arrow $\alpha$ ending at $v$.
 By construction, $N$ is a separated monic representation
 and $N_i$ is projective for each $i\in Q_0$.
 Then $N$ is projective in $\Gproj \Lambda Q$ by Corollary \ref{pro}.
\medskip

Since $\sigma_X:\Omega_X\to C_X $ is a monomorphism and $M_v$ is an injective $\Lambda$-module,
then there exists a map $u: C_X\to M_v$ such that $h_v=u\sigma_X$.
Since $h: X(v)\to M$ is a homomorphism, for each $\alpha$ ending at $v$,
$ h_{s(\alpha)} \left(\begin{smallmatrix} k _{X,\alpha^{op}} \\ \sigma_X\\0\end{smallmatrix}\right)
  = M_{\alpha^{op}}  h_v = M_{\alpha^{op}}  u\sigma_X$,
  that is, $h^X_{s(\alpha)}k_{X,\alpha^{op}}= (M_{\alpha^{op}}  u-h_{\alpha^{op}})\sigma_X$.

  \medskip
  Now we use the form of $X(v)$ as described in Subsection~\ref{sec-reflection-objects}:
  Hence $(h^X_{s(\alpha)})_\alpha k_X=(M_{\alpha^{op}}  u-h_{\alpha^{op}})_\alpha\sigma_X$
  where $\alpha$ runs over all arrows ending at $v$ in $Q_1$.
  Since $X_v$ is a pushout of  $(k_X, \sigma_X)$,
  there is an unique map $(z_{\alpha^{op}})_{\alpha }: X_v\to \bigoplus\limits_{e(\alpha)= v}M_{s(\alpha)}$  such that
  $$(h^X_{s(\alpha)})_{\alpha }=(z_{\alpha^{op}})_{\alpha } u_X, (M_{\alpha^{op}}  u-h_{\alpha^{op}})_\alpha= (z_{\alpha^{op}})_{\alpha } t_X.\ \ \ \ \ \ (3.1)$$
  By the equations $f(v)_v=\Omega_f=g_vh_v,\  \ \sigma_Y\Omega_f=C_f\sigma_X$ and $h_v=u\sigma_X $,
  we have $(\sigma_Y g_vu-C_f) \sigma_X=0,$
  hence there is a map $\gamma: X_{\bar v}\to C_Y $ such that $\sigma_Y g_vu-C_f=\gamma\delta_X$.
  Since $g$ is a homomorphism, $f(v)=gh$ and the equation $(3.1)$, we have that $$u_Y(g^Y_{s(\alpha)})_\alpha(z_{\alpha^{op}})_{\alpha} t_X=t_Y\gamma\delta_X+t_YC_f-u_Y\omega_f=t_Y\gamma\delta_X+f_vt_X.$$
  Hence
  $$[u_Y(g^Y_{s(\alpha)})_\alpha(z_{\alpha^{op}})_{\alpha}- t_Y\gamma\pi_X]u_X=u_Y(f_{s(\alpha)})_\alpha$$
  $$[u_Y(g^Y_{s(\alpha)})_\alpha(z_{\alpha^{op}})_{\alpha}- t_Y\gamma\pi_X] t_X=t_YC_f-u_Y\omega_f.$$
  Since $X_v$ is the pushout of $(k_X, \sigma_X)$,
$ u_Y (f_{s(\alpha)})_{\alpha}k_X=(t_YC_f-u_Y\omega_f)\sigma_X$,
  and $ f_vu_X=u_Y(f_{s(\alpha)})_{\alpha}, f_vt_X=t_YC_f-u_Y\omega_f$,
  then by the uniqueness of $f_v$, we have
  $f_v=u_Y(g^Y_{s(\alpha)})_\alpha(z_{\alpha^{op}})_{\alpha}- t_Y\gamma\pi_X$.
  Let
  $h'_v=\left(\begin{smallmatrix}(z_{\alpha^{op}})_\alpha\\ \gamma \pi_X\end{smallmatrix}\right): X_v\to N_v,\    g'_v=\left(\begin{smallmatrix}u_Y(g^Y_{s(\alpha)})_\alpha,& -t_Y \end{smallmatrix}\right):N_v\to Y_v $, then $f_v=g'_vh'_v$.
    Let $g'=(g'_i)$
    and $h'=(h'_i)$, where
    $$g'_i = \begin{cases} g'_v, &  \ i =v; \\
      g^Y_{i },  &    i \neq v.\end{cases},\ h'_i = \begin{cases} h'_v, &  \ i =v; \\
      h^X_{i},  &    i \neq v.\end{cases}, $$
    then $g'$ and $h'$ are homomorphisms and $f=g'h'$ which means that $f$ factorizes through the projective object $N$ in $ \Gproj  \Lambda Q $.

\vskip5pt
Claim $2$: $F(v) $ is full.
\vskip5pt
In fact, let $g: X(v)\to Y(v)$ be a homomorphism in $ \underline{ \Gproj} \Lambda Q(v)$, we need to find a homomorphism  $f: X \to Y $  in $   \Gproj \Lambda Q $, such that the class of $f(v)$ is equal to $ g$ in $ \underline{ \Gproj} \Lambda Q(v)$.
For each arrow $\beta:v\to i$ in $Q(v)$, we write
$ g_{i}: X(v)_{i}\to Y(v)_{i}$ as
$$\left(\begin{smallmatrix}
g_{iXY}&g_{i\beta Y}& (g_{ip Y})_p\\
g_{iX\beta}&g_{i\beta \beta }& (g_{ip\beta })_p\\
(g_{iXq})_q&(g_{i\beta q })_q& (g_{ipq})_{pq}\\
\end{smallmatrix}\right): X_{i}\oplus C_{X, \beta}\oplus\bigoplus\limits_{p\neq \beta}C_{X,p} \to Y_{i}\oplus C_{Y, \beta}\oplus\bigoplus\limits_{q\neq \beta}C_{Y,q}$$
where $p$ and $q$ run over all non-trival paths from $v$ to $i$ excluding $\beta$.
Since $g_{i}X(v)_\beta=Y(v)_\beta g_v$ for each arrow $\beta:v\to i$ in $Q(v)$, we have $  g_{iXY}k_{X,\beta} +g_{i\beta Y} \sigma_X=k_{Y,\beta} g_v,  g_{iX\beta}k_{X,\beta} +g_{i\beta\beta} \sigma_X=\sigma_Y g_v,  (g_{iXq})_q k_{X,\beta}+ (g_{i\beta q })_q\sigma_X=0.$  So
$$( g_{e(\beta)XY})_\beta k_X +(g_{t(\beta)\beta Y} )_\beta \left(\begin{smallmatrix}1\\  \vdots\\ 1\end{smallmatrix}\right)  \sigma_X=k_Y g_v,$$  $$( g_{t(\beta)X\beta})_\beta k_X+(g_{t(\beta)\beta\beta})_\beta \left(\begin{smallmatrix}1\\  \vdots\\ 1\end{smallmatrix}\right)  \sigma_X= \left(\begin{smallmatrix}1\\  \vdots\\ 1\end{smallmatrix}\right)  \sigma_Y g_v,\ \ \ \ \ \ \ \ (3.2)$$ where $\beta$ runs over all arrows starting at $v$.
Since $u_Y k_Y=t_Y\sigma_Y$, we have
$$ [\frac{1}{m}\left(\begin{smallmatrix} t_Y,&  \cdots,&  t_Y\end{smallmatrix}\right)( g_{t(\beta)X\beta})_\beta  -u_Y( g_{t(\beta)XY})_\beta ]k_X $$ $$=[u_Y(g_{t(\beta)\beta Y} )_\beta \left(\begin{smallmatrix}1\\  \vdots\\ 1\end{smallmatrix}\right) -\frac{1}{m}\left(\begin{smallmatrix} t_Y,&  \cdots,&  t_Y\end{smallmatrix}\right) (g_{t(\beta)\beta\beta} )_\beta \left(\begin{smallmatrix}1\\  \vdots\\ 1\end{smallmatrix}\right)  ]  \sigma_X$$
where $m$ is the number of arrows starting at $v$.
Since $X_v$ is the pushout of $(k_X, \sigma_X)$, then there is a unique map $s: X_v\to Y_v$ such that
$$s u_X=\frac{1}{m}\left(\begin{smallmatrix} t_Y,&  \cdots,&  t_Y\end{smallmatrix}\right)( g_{t(\beta)X\beta})_\beta  -u_Y( g_{t(\beta)XY})_\beta, $$  $$st_X= u_Y(g_{t(\beta)\beta Y} )_\beta \left(\begin{smallmatrix}1\\  \vdots\\ 1\end{smallmatrix}\right) -\frac{1}{m}\left(\begin{smallmatrix} t_Y,&  \cdots,&  t_Y\end{smallmatrix}\right) (g_{t(\beta)\beta\beta} )_\beta \left(\begin{smallmatrix}1\\  \vdots\\ 1\end{smallmatrix}\right). $$
Since $u_X$ is a monomorphism and $C_{Y,\beta}=C_Y$ is an injective $\Lambda$-module, then there is a $\Lambda$-map $\gamma=(\gamma_\beta)_\beta: X_v\to \bigoplus\limits_{\beta}C_{Y,\beta}$ where $\beta $ runs over all arrows starting at $v$, such that $( g_{t(\beta)X\beta})_\beta =\gamma u_X. $
Let $$f_i = \begin{cases}g_{iXY}, &  \ i \neq v; \\
 \frac{1}{m}(t_Y,\cdots, t_Y)\gamma-s,  &    i=v,\end{cases}$$
 then $f_vu_X=u_Y( g_{t(\beta)XY})_\beta.$ Moreover, $f=(f_i):X\to Y$ is a homomorphism following from the fact that $g:X(v)\to Y(v)$ is a homomorphism and  $f_vu_X=u_Y( g_{t(\beta)XY})_\beta$.
By the construction of $f(v)$, we know that $f(v)_v=\Omega_f$ and$$ f(v)_i=\left(\begin{smallmatrix} f_i& (Y_{p'}\omega_{f, \beta})_{p } \\ 0&C_fE \end{smallmatrix}\right)$$ for each $i\neq v$ where $p$ runs over all non-trival paths from $v$ to $i$ and $\beta$ is the starting arrow of $p$ such that $p= p'\beta$  and the size of $E$ is determined by the number of nontrival paths from $v$ to $i$.

 \medskip

Sub-claim: the class of $f(v)$ is equal to $g$ in $ \underline{ \Gproj} \Lambda Q(v)$. That  is to say, there is  a projective object $P$ and two homomorphisms $h: X(v)\to P$ and $k: P\to Y(v)$ such that $g-f(v)=kh$ in $\Gproj  \Lambda Q(v)$.

 \medskip
  Let $T=\frac{1}{m}\left(\begin{smallmatrix}1,&  \cdots,& 1\end{smallmatrix}\right)\gamma t_X+\frac{1}{m}\left(\begin{smallmatrix} 1,&  \cdots,& 1\end{smallmatrix}\right) (g_{t(\beta)\beta\beta} )_\beta \left(\begin{smallmatrix}1\\  \vdots\\ 1\end{smallmatrix}\right)$, where $\beta$ runs over all arrows starting at $v$.
Since $\delta_Y=\pi_Y t_Y$, $ \pi_Y u_Y=0$ and $u_Y\omega_f=t_YC_f-f_vt_X$, then we have $\delta_Y(C_f-T)=0$.
Hence there is a homomorphism $l: C_X\to \Omega_Y$ such that  $\sigma_Y l=C_f-T. $
 Since $C_f\sigma_X=\sigma_Y\Omega_f, t_X\sigma_X=u_Xk_X$, $\gamma u_X=(g_{t(\beta)X\beta})_\beta$ and the equation $(3.2)$, then there holds $\sigma_Yl \sigma_X=\sigma_Y(\Omega_f-g_v) $. So $ l \sigma_X= \Omega_f-g_v$ follows from the fact that $\sigma_Y$ is a monomorphism.
\medskip

Let $P=P_v(C_X)\oplus \bigoplus\limits_{s(\beta)=v }P_{t(\beta)}(C_Y), h=(h_i)$ and $k=(k_i)$, where $ h_v=\sigma_X, k_v=-l$ and for each $i\neq v$,
$$ h_i=\left(\begin{smallmatrix} 0&E\\ (g_{iXq})_q& (g_{ipq})_{pq}-TE\end{smallmatrix}\right), k_i=\left(\begin{smallmatrix} (g_{ipY}-Y_{p'}\omega_{f,\beta})_p&0\\ (T-C_f)E& E\end{smallmatrix}\right),   $$ where $p$ and $q$ runs over all non-trival paths from $v$ to $i$, and $p'$ is the subpath of $p$ such that $p=p'\beta$ with $\beta\in Q(v)_1$, $E$ is the identity matrix whose size is determined by the number of non-trival paths from $v$ to $i$. Hence $k_vh_v=-l\sigma_X=g_v-\Omega_f=g_v-f(v)_v.$ For each $i\neq v$, by the multiplication of matrix, we have $k_ih_i=g_i-f(v)_i.$ So for each $i\in Q(v)_0,$ $k_ih_i=g_i-f(v)_i.$ Moreover, we can prove that $k$ and $h$ are homomorphisms in $\Gproj \Lambda Q(v).$ Hence $F$ is full.

\end{proof}

  Iterated application of the previous result yields

\begin{thm}
  Suppose that $\Lambda$ is a selfinjective algebra, $Q$ and $\widetilde Q$ are finite acyclic quivers, and
  that the quiver $\widetilde Q$ can be obtained from the quiver $Q$ by a sequence
  of reflections.  Then $$\uGproj \Lambda Q\;\simeq\; \uGproj\Lambda\widetilde Q.$$
\end{thm}

\vskip10pt

  %--------------------------------------------------------------------------------
  \section{The Quiver $\mathbb A_3$ revisited}
  %--------------------------------------------------------------------------------
  \label{section-reflections}

  In this section we aim to show the following statement.

  \begin{prop}
    \label{proposition-twelve}
    Let $\Lambda$ be a selfinjective algebra and $Q$ be the quiver $Q:1\to2\to3$.
    Suppose that $M$ is a representation for $Q$ with coefficients in $\Lambda$
    such that none of the branches $M_i$, $i\in Q_0$, of $M$
    has a nonzero injective direct summand.

    Let $A=\Mimo(M)$ be the corresponding Gorenstein-projective representation
    and let $B=A(1,2,3,1,2,3,1,2,3,1,2,3)$ be obtained from $A$ by performing
    twelve reflections.  Then
    $$B\cong\Omega_\Lambda^2A.$$
  \end{prop}

   The functor $\Omega_\Lambda$ is the inverse of the suspension functor $[1]$
   for the stable category $\umod\Lambda$. The functor $\Omega_\Lambda^2$
   is applied pointwise on the $\Lambda Q$-module $A$.

  \begin{cor}
    Suppose that $\Lambda$ is a selfinjective algebra such that
    the functor $\Omega_\Lambda^2$ is naturally equivalent to the identity functor
    on the stable category, and $Q$ is  the quiver $Q:1\to2\to3$.
    If $M$, $A$, $B$ are as in Proposition~\ref{proposition-twelve}, then
    $$B\cong A.$$
  \end{cor}

  \subsection{
    On the geometry of the octahedron in the octahedral axiom}
  %--------------------------------------------------------------------------------

  \label{sec-octa}
    
  Consider a pair $(u,v)$ of composable maps in a triangulated category.  The three maps
  $u, v, v\circ u$ give rise to three exact triangles.  The octahedral axiom states that
  the mapping cones of these three triangles
  are the vertices of a fourth exact triangle such that ``everything commutes''.

  \medskip
  The four exact triangles form four of the eight faces of an octahedron.
  The other four faces consist of pairs of
  composable maps which together with their composition form a non-oriented triangle.
  These four triangles are arranged such that the four products form an
  oriented cycle.
  
  \smallskip
  Note that an octahedron contains three squares such that any two are perpendicular to each other.
  One square is given by the oriented cycle of the four products just mentioned.
  The other two squares are perpendicular to it, and to each other, and are
  required by the axiom to be commutative squares.
  The starting vertex of one square is the end vertex of the other,
  and conversely.  Those two vertices, $Y'$ and $Y$ in the picture below,
  define an axis which is perpendicular to the oriented cycle.

  \smallskip
  Consider the oriented cycle as the ``equator'' of the octahedron.
  We remark that the orientation given
  by the oriented cycle does not distinguish a ``North pole'' over a ``South pole'' (or conversely)
  as the axiom is symmetric
  with respect to reflection on the plane given by the equator.  But the orientation of the oriented
  cycle does define a cyclic ordering of the four exact triangles involved in the octahedron.

  \begin{obs}
    \label{obs-symmetry}
    Consider the octahedron in the octahedral axiom as an oriented graph $G$.
    The symmetry group of $G$ is cyclic of order four.
    It is generated by a rotation-reflection $\rho$ given by rotating the octahedron about
    the axis mentioned above by $90^\circ$ in the direction of the orientation of the cyclically oriented square,
    and then reflecting the octahedron about the plane given by that square.
  \end{obs}

%  \smallskip
%  We will see below that any three subsequent reflections effect a rotation of the octahedron
%  about its axis by $90^\circ$ in the direction given by the oriented cycle.

  \smallskip
  We present in Figure~\ref{figure-octa} the octahedron using a cylindrical projection
  from the axis which is perpendicular to the cyclically oriented square.
  Thus, the equator given by the square consists of
  the four products and is given by the horizontal line in the middle.
  Above and below the equator, exact triangles and triangles given by a pair of composable maps alternate.
  Under the cylindrical projection,
  the generator $\rho$ in Observation~\ref{obs-symmetry} appears as a glide-reflection.
  Note that while $\rho^4$ acts as the identity map on the octahedron, it operates as the square $[2]$ of the
  suspension on its exact triangles.  (This is to be expected since each pole
  $Y$ or $Y'$ occurs both as starting term
  and as end term of two exact triangles.)

  \smallskip
  In the diagram, vertices and arrows are labelled as in \cite[Page 3]{Hap1}.
  The four faces of the octahedron which are given by pairs of
  composable maps occur as rectangles which contain entries;
  the other four faces are exact triangles which appear as empty rectangles.

%  \smallskip
%  We project the octahedron against a cylinder around this axis.
  \begin{figure}[ht]
  \begin{tikzcd}[column sep=scriptsize, row sep=large]
    & Y'[-1] \arrow[d, swap, "{k'[-1]}"]
    & Y'
    & Y' \arrow[d, swap, "g"]
    & Y'[1]
    & Y'[1] \arrow[d, swap, "{k'[1]}"]
    & Y'[2]
    & \\
    \arrow[r]
    & X \arrow[d, "u"] \arrow[r, swap, "v\circ u"]
    & Z \arrow[u, swap, "k"] \arrow[r, "g\circ k"]
    & X' \arrow[d, "j'"] \arrow[r, swap, "{i[1]\circ j'}"]
    & Z'[1] \arrow[u, swap, "{f[1]}"] \arrow[r,"{(k'\circ f)[1]}"]
    & X[2] \arrow[d,"{u[2]}"] \arrow[r, swap, "{(v\circ u)[2]}"]
    & Z[2] \arrow[u, swap, "{k[2]}"] \arrow[r]
    & \phantom{.} \\
    & Y
    & Y \arrow[u,"v"]
    & Y[1]
    & Y[1] \arrow[u,"{i[1]}"]
    & Y[2]
    & Y[2] \arrow[u,"{v[2]}"]
    &
  \end{tikzcd}
  \caption{The octahedron of the axiom in cylindrical projection}
  \label{figure-octa}
  \end{figure}

  %--------------------------------------------------------------------------------
  \subsection{The first reflection}
  %--------------------------------------------------------------------------------

  Let $Q$ be the quiver of type $\mathbb A_3$ in linear orientation,
  $$Q:\qquad1\;\stackrel\alpha\longrightarrow \;2\;\stackrel\beta\longrightarrow\;3$$
  and $\Lambda$ a selfinjective algebra.

  \medskip
  Consider $M$, a representation of $Q$  over $\Lambda$, such that no component $M_i$
  has a non-zero injective direct summand. Thus, $M$ consists of a pair
  of composable maps in the stable category $\umod\Lambda$:
  $$M: \qquad X\;\stackrel u\longrightarrow\;Y\stackrel v\longrightarrow\;Z$$
  We fix the octahedron given by this pair of composable maps and call it $\mathcal O$.

  \medskip
  Then $A=\Mimo(M):\;A_1\stackrel\alpha\to A_2\stackrel\beta\to A_3$
  is an object in $\Gproj\Lambda Q$, we apply reflections to $A$.

  \medskip
  We first reflect at 3.  Recall that $A(3)=\Mimo(A'(3))$ where $A'(3)$ is a
  representation of $Q(3):\;1\stackrel\alpha\to 2\stackrel\gamma\leftarrow3$ given as
  $$A'(3):\qquad A_1\stackrel\alpha\longrightarrow A_2\stackrel k\longleftarrow\Omega_A$$
  where the map
  $k$ representing $\gamma$ is as in the following diagram.
  \begin{center}\begin{tikzcd}
      0 \arrow[r] & \Omega_A \arrow[r, "\sigma_A"] \arrow[d,swap,"k"]
      & C_A \arrow[d, "t_A"] \arrow[r,"\delta_A"]
      & A_{\bar 3}\arrow[r] \arrow[d, equal] & 0 \\
      0 \arrow[r] & A_2 \arrow[r, "\beta"] & A_3 \arrow[r,"\pi"] & A_{\bar 3} \arrow[r] & 0
  \end{tikzcd}\end{center}
  We assume that $\Omega_A$ has no nonzero injective direct summand.
  In the stable category, the diagram represents the triangle
  $$A_2\stackrel\beta\longrightarrow A_3\stackrel \pi\longrightarrow A_{\bar 3}
  \stackrel{-k[1]}\longrightarrow$$

  which is equivalent to the triangle in the octahedron $\mathcal O$,
  $$Y\stackrel v\longrightarrow Z\stackrel {g\circ k}\longrightarrow
  X'\stackrel{j'}\longrightarrow,$$
  since the maps $A_2\to Y$, $A_3\to Z$ are isomorphisms in $\umod\Lambda$ which make
  the left square commutative; hence $\Omega_A\cong X'[-1]$ in $\umod\Lambda$
  and the square is commutative.
  \begin{center}\begin{tikzcd}
      \Omega_A \arrow[r,"k"] \arrow[d,swap,"\cong"] & A_2 \arrow[d,"\cong"] \\
      X'[-1] \arrow[r,swap,"{j'[-1]}"] & Y
  \end{tikzcd}\end{center}

  We obtain a representation of $Q(3)$ consisting of objects and morphisms
  in the octahedron $\mathcal O$.
  $$M(3):\qquad X\;\stackrel u\longrightarrow \;Y\;\stackrel {j'[-1]}\longleftarrow\;X'[-1].$$

  %--------------------------------------------------------------------------------
  \subsection{The second reflection}
  %--------------------------------------------------------------------------------

  Next, we reflect at vertex 2.  The representation for $Q(3)$,
  $$U=A(3)=\Mimo(A'(3)):\qquad U_1\;\stackrel\alpha\longrightarrow\; U_2\;
  \stackrel\gamma\longleftarrow \;U_3$$
  is separated monic, so we can compute $U'(2)$ using the following diagram.
  We assume that $\Omega_U$ has no nonzero injective direct summand.
  \begin{center}\begin{tikzcd}
      0 \arrow[r] & \Omega_U \arrow[r, "\sigma_U"] \arrow[d,swap,"{k_1\choose k_3}"]
      & C_U \arrow[d, "t_U"] \arrow[r,"\delta_U"]
      & U_{\bar 2}\arrow[r] \arrow[d, equal] & 0 \\
      0 \arrow[r] & U_1\oplus U_3 \arrow[r, "{(\alpha,\gamma)}"]
      & U_2 \arrow[r,"\pi"] & U_{\bar 2} \arrow[r] & 0
  \end{tikzcd}\end{center}
  Hence, $U'(2)$ is the representation for the quiver $Q(2,3):\;1\stackrel{\alpha^\op}\leftarrow
  2\stackrel\beta\to 3$ given by
  $$U'(2):\qquad U_1\;\stackrel{k_1}\longleftarrow\;\Omega_U\;\stackrel{k_3}\longrightarrow\;U_3.$$

  \medskip
  The interesting step is to detect the maps $k_1:\Omega_U\to U_1$ and $k_3:\Omega_U\to U_3$
  in the octahedron.  We first deal with $k_3$.
  Using that $U$ is separated monic, we obtain the commutative diagram with
  exact rows and columns.
  \begin{center}\begin{tikzcd}
      & 0 \arrow[d] & 0 \arrow[d] & & \\
      & U_1 \arrow[r,equal] \arrow[d,"{{1 \choose 0}}"] & U_1 \arrow[d,"\alpha"] & & \\
      0 \arrow[r] & U_1\oplus U_3 \arrow[r, "{(\alpha,\gamma)}"] \arrow[d,swap,"{(0,1)}"]
      & U_2 \arrow[d, "{{\rm can}_{U_1}}"] \arrow[r,"\pi"]
      & \frac{U_2}{U_1\oplus U_3} \arrow[r] \arrow[d, equal] & 0 \\
      0 \arrow[r] & U_3 \arrow[r, "{{\rm can}_{U_1}\circ\gamma}"] \arrow[d]
      & \frac{U_2}{U_1} \arrow[r,"{\rm can}"] \arrow[d] & \frac{U_2}{U_1\oplus U_3} \arrow[r] & 0 \\
      & 0 & 0 & &
  \end{tikzcd}\end{center}
  Note that the middle row of this diagram is the bottom row of the diagram
  above.  Combining the two diagrams yields the commutative diagram with exact rows.
  \begin{center}\begin{tikzcd}
      0 \arrow[r] & \Omega_U \arrow[r, "\sigma_U"] \arrow[d,swap,"{k_1\choose k_3}"]
      & C_U \arrow[d, "t_U"] \arrow[r,"\delta_U"]
      & U_{\bar 2}\arrow[r] \arrow[d, equal] & 0 \\
      0 \arrow[r] & U_1\oplus U_3 \arrow[r, "{(\alpha,\gamma)}"] \arrow[d,swap,"{(0,1)}"]
      & U_2 \arrow[d, "{{\rm can}_{U_1}}"] \arrow[r,"\pi"]
      & U_{\bar 2} \arrow[r] \arrow[d, equal] & 0 \\
      0 \arrow[r] & U_3 \arrow[r, "{{\rm can}_{U_1}\circ\gamma}"]
      & \frac{U_2}{U_1} \arrow[r,"{\rm can}"] & U_{\bar 2} \arrow[r] & 0
  \end{tikzcd}\end{center}
  By composing the vertical maps and omitting the middle row we obtain the following diagram.
  \begin{center}\begin{tikzcd}
      0 \arrow[r] & \Omega_U \arrow[r, "\sigma_U"] \arrow[d,swap,"k_3"]
      & C_U \arrow[d, "t'_U"] \arrow[r,"\delta_U"]
      & U_{\bar 2}\arrow[r] \arrow[d, equal] & 0 \\
      0 \arrow[r] & U_3 \arrow[r, "{{\rm can}_{U_1}\circ\gamma}"]
      & \frac{U_2}{U_1} \arrow[r,"{\rm can}"] & U_{\bar 2} \arrow[r] & 0
  \end{tikzcd}\end{center}
  Since $C_U$ is a projective cover, this diagram represents the triangle in $\umod\Lambda$:
  $$U_3\;\stackrel{{\rm can}_{U_1}\circ\gamma}\longrightarrow \;\frac{U_2}{U_1}\;
  \stackrel{\rm can}\longrightarrow\;U_{\bar 2}\stackrel{-k_3[1]}\longrightarrow$$
  Consider the first morphism.  The map $\gamma:U_3\to U_2$ is equivalent to
  $j'[-1]:X'[-1]\to Y$ (compare $A(3)$ and $M(3)$ above);
  the canonical map ${\rm can}_{U_1}:U_2\to U_2/U_1$,
  which is the cokernel of $\alpha$, is equivalent to the map $i$ in the octahedron
  since $\alpha$ is equivalent to $u$ and
  $X\stackrel u\to Y\stackrel i\to Z'\to$ is a triangle in which $i$ follows $u$.
  Hence the above triangle is equivalent to the following triangle in the octahedron:
  $$X'[-1]\;\stackrel{-i\circ j'[-1]}\longrightarrow Z'\stackrel f\longrightarrow Y'
  \stackrel g\longrightarrow$$
  Hence the map $k_3:\Omega_U\to U_3$ is equivalent to
  the map $g[-1]:Y'[-1]\to X'[-1]$.

  \medskip
  In a similar way we obtain $k_1$.  Consider the first diagram in this subsection
  and the second diagram with the roles of $U_1$ and $U_3$ exchanged.
  Again, the bottom row of the first diagram is just the middle row of the second,
  so we obtain the following diagram by
  composing the vertical maps and omitting the middle row.
  \begin{center}\begin{tikzcd}
      0 \arrow[r] & \Omega_U \arrow[r, "\sigma_U"] \arrow[d,swap,"k_1"]
      & C_U \arrow[d, "t''_U"] \arrow[r,"\delta_U"]
      & U_{\bar 2}\arrow[r] \arrow[d, equal] & 0 \\
      0 \arrow[r] & U_1 \arrow[r, "{{\rm can}_{U_3}\circ\alpha}"]
      & \frac{U_2}{U_3} \arrow[r,"{\rm can}"] & U_{\bar 2} \arrow[r] & 0
  \end{tikzcd}\end{center}
  We read off that
  $$U_1\;\stackrel{{\rm can}_{U_3}\circ \alpha}\longrightarrow\;\frac{U_2}{U_3}\;
  \stackrel{{\rm can}}\longrightarrow\;U_{\bar 2}\;\stackrel{k_1[1]}\longrightarrow$$
  is a triangle.  By comparing first maps, we see that it is equivalent to the
  following triangle in the octahedron.
  $$X\;\stackrel{ v\circ u}\longrightarrow \;Z\; \stackrel k\longrightarrow
  \;Y'\;\stackrel{k'}\longrightarrow$$
  Hence $k_1$ is equivalent to the map $k'[-1]$ in the octahedron.

  \medskip
  In conclusion, the representation $U'(2)$, and hence $A(2,3)=\Mimo(U'(2))$,
  is equivalent in $\umod\Lambda$ to the following representation which consists
  of objects and morphisms in $\mathcal O$.
  $$M(2,3):\qquad X\; \stackrel{k'[-1]}\longleftarrow \;Y'[-1] \;
  \stackrel{g[-1]}\longrightarrow \; X'[-1].$$
  Note that, together, the maps $u$, $j'[-1]$ defining $M(3)$
  and the maps $k'[-1]$ and $g[-1]$ defining $M(2,3)$ form a square in the octahedron
  $\mathcal O$
  \begin{center}\begin{tikzcd}[row sep=small]
      & Y'[-1] \arrow[dl,swap, bend right, "{k'[-1]}"] \arrow[dr, bend left, "{g[-1]}"] & \\
      X \arrow[dr, bend right, swap,"u"] & & X'[-1] \arrow[dl, bend left,"{j'[-1]}"] \\
      & Y &
  \end{tikzcd}\end{center}
  the commutativity of which is required by the axiom.

  %--------------------------------------------------------------------------------
  \subsection{The third reflection}
  %--------------------------------------------------------------------------------

  In order to return to a representation of the quiver $Q$, it remains to reflect at vertex 1.
  We are given a separated monic representation for the quiver $Q(2,3)$:
  $$V=A(2,3):\qquad V_1\;\stackrel{\delta}\longleftarrow\; V_2
  \;\stackrel\beta\longrightarrow\; V_3$$
  We may assume that $V_2$ has no nonzero injective direct summand.
  For the reflection, we consider the diagram in which $\Omega_V$, again,
  is supposed to have no nonzero injective direct summand.
  \begin{center}\begin{tikzcd}
      0 \arrow[r] & \Omega_V \arrow[r, "\sigma_V"] \arrow[d,swap,"\ell"]
      & C_V \arrow[d, "t_V"] \arrow[r,"\delta_V"]
      & V_{\bar 1}\arrow[r] \arrow[d, equal] & 0 \\
      0 \arrow[r] & V_2 \arrow[r, "\delta"] & V_1 \arrow[r,"\pi"] & V_{\bar 1} \arrow[r] & 0
  \end{tikzcd}\end{center}
  We obtain the representation $V'(1)$ for the quiver $Q$:
  $$V'(1):\qquad \Omega_V\;\stackrel \ell\longrightarrow \;V_2\;
  \stackrel \beta\longrightarrow \;V_3$$
  We locate a map equivalent to $\ell$ in the octahedron.  Recall that $\delta$
  is equivalent to the corresponding map $k'[-1]:Y'[-1]\to X$ in $M(2,3)$
  which occurs in the triangle
  $$Z[-1]\;\stackrel{k[-1]}\longrightarrow\;Y'[-1]\;\stackrel{k'[-1]}\longrightarrow
  X \;\stackrel{v\circ u}\longrightarrow .$$

  Here is the representation $M(1,2,3)$ of $Q$ with objects and morphisms
  in the octahedron $\mathcal O$.  It is equivalent in $\umod\Lambda$ to
  $A(1,2,3)=\Mimo(V'(1))$.
  $$ M(1,2,3): \qquad Z[-1]\;\stackrel{k[-1]}\longrightarrow \;Y'[-1]\;
  \stackrel {g[-1]}\longrightarrow \;X'[-1]$$

  %--------------------------------------------------------------------------------
  \subsection{Summary of the first three reflections}
  %--------------------------------------------------------------------------------

  In the octahedron, $M(1,2,3)$ is represented by a pair of composable maps.
  We compare this pair $(k[-1],g[-1])$ to the pair $(u,v)$ for $M$.
  Up to a shift by $-1$, the pair $(k[-1],g[-1])$ is adjacent to the original pair $(u,v)$,
  but on the opposite side of the equator of the octahedron.
  Hence up to the shift $[-1]$, three subsequent reflections are represented by the
  rotation-reflection $\rho$ of the octahedron, see Observation~\ref{obs-symmetry}.
  
  %--------------------------------------------------------------------------------
  \subsection{Twelve reflections}
  %--------------------------------------------------------------------------------

  \sloppypar
  We continue this process to determine the representation obtained from $M$
  after 12 reflections, this is $M(1,2,3,1,2,3,1,2,3,1,2,3)$.
  For this we can repeat the three above steps four times
    and trace the effect on the octahedron in Figure~\ref{figure-octa}.
    Recall that $M$ is given by a pair of composable maps,
    $$M: \qquad X\stackrel u\longrightarrow Y\stackrel v\longrightarrow Z.$$
    We have seen that $M(1,2,3)$ is given by the following pair of composable maps,
    $$M(1,2,3):\qquad Z[-1]\stackrel{k[-1]}\longrightarrow Y'[-1]\stackrel{g[-1]}\longrightarrow X'[-1].$$
    Repeating this process shows that $M(1,2,3,1,2,3)$ is as follows.
    $$M(1,2,3,1,2,3):\qquad X'[-2]\stackrel{j'[-2]}\longrightarrow Y[-1]\stackrel{i[-1]}\longrightarrow Z'[-1]$$
    Similarly,
    $$M(1,2,3,1,2,3,1,2,3):\qquad Z'[-2]\stackrel{f[-2]}\longrightarrow
    Y'[-2]\stackrel{k'[-2]}\longrightarrow X[-1].$$
    Finally, we obtain
    $$M(1,2,3,1,2,3,1,2,3,1,2,3):\qquad X[-2]\stackrel{u[-2]}\longrightarrow
    Y[-2]\stackrel{v[-2]}\longrightarrow Z[-2].$$
    In conclusion, $M(1,2,3,1,2,3,1,2,3,1,2,3)\cong M[-2]$ is obtained by applying $\rho$ four times
    to the composable pair $(u,v)$ (resulting in $(u[2],v[2])$) and by applying the inverse suspension $[-1]$
    also four times.  This finishes
    the proof of Proposition~\ref{proposition-twelve}.
    We illustrate the construction in Section~\ref{sec-12} in an example.

%--------------------------------------------------------------------------------
\section{Examples}
%--------------------------------------------------------------------------------
\label{section-examples}

  The examples in this section are Gorenstein-projective representations
  for $\Lambda Q$ where $Q$ is a quiver of type $\mathbb A_3$ in various
  orientations and $\Lambda=k[T]/(T^n)$ the bounded polynomial ring
  where $n$ is $2$ or $3$.

  \subsection{Some Gorenstein-projective representations}

\medskip
In the case where $n=2$ we present a table to show
that the category of Gorenstein-projective
quiver representations may behave nicer than either the
monic representations, or all representations of the quiver
in the sense that the Gorenstein-projective representations form a
Frobenius category, and the stable part of this category does not
depend on the orientation of the quiver, as we have seen in
Theorem~\ref{theorem-intro1}.

\medskip
For $Q$ a quiver of type $\mathbb A_3$ in any orientation
and $\Lambda=k[\varepsilon]=k[T]/(T^2)$,
the category of $\Lambda Q$-modules has finite representation type.
We list in Table~\ref{table1} for each orientation:
the number of indecomposables,
the shape of the stable part of the AR-quiver, the number of monic
indecomposables, the number of Gorenstein-projective indecomposables,
and the shape of the stable part of the AR-quiver for $\Gproj\Lambda Q$.

\medskip
\begin{center}
  \begin{table}
    \label{table1}
\begin{tabular}{|c|c|c|c|c|c|} \hline
  $\rule[-1mm]{0mm}{6mm}Q$ & $\#\ind \Lambda Q$ & $\overline{\Gamma_{\Lambda Q}}$
  & $\#\mon\Lambda Q$ & $\#\Gp\Lambda Q$
  & $\overline{\Gamma_{\Gp\Lambda Q}}$ \\ \hline\hline
  $\rule[-2mm]{0mm}{9mm}\bullet\to\bullet\to\bullet$ & 36
  & $\D\frac{\mathbb Z\mathbb D_6}{\tau^5}$
  & 9 & 9 & $\D\frac{\mathbb Z\mathbb A_3}{\tau^2\varphi}$ \\ \hline
  $\rule[-2mm]{0mm}{9mm}\bullet\to\bullet\gets\bullet$ & 42
  & $\D\frac{\mathbb Z\mathbb E_6}{\tau^6\varphi}$
  & 16 & 9 & $\D\frac{\mathbb Z\mathbb A_3}{\tau^2\varphi}$\\ \hline
  $\rule[-2mm]{0mm}{9mm}\bullet\gets\bullet\to\bullet$ & 42
  & $\D\frac{\mathbb Z\mathbb E_6}{\tau^6\varphi}$
  & 9 & 9 & $\D\frac{\mathbb Z\mathbb A_3}{\tau^2\varphi}$\\ \hline
\end{tabular}\\[0.5ex]
Table~\ref{table1}: Representations of $Q$ over the ring $\Lambda=k[\varepsilon]=k[T]/(T^2) $
\end{table}
\end{center}
We picture the Auslander-Reiten quivers for the three categories of
Gorenstein-projective representations.  In each case we describe in detail
how the objects give rise to the icons which occur
in the Auslander-Reiten quiver.
The Auslander-Reiten quivers are obtained using universal coverings
as in \cite[Part A]{RS2}.

\medskip
Recall that a nilpotent
linear operator $V$ consists of a finite dimensional $k$-vector
space, also denoted $V$,
together with a linear map $T:V\to V$ such that $T^n=0$ for some
$n\in\mathbb N$.  Up to isomorphy, $V$ is given uniquely by a partition
recording the dimensions of the indecomposable direct summands,
that is, the sizes of the Jordan blocks of $T$.
We picture the parts of the partition
as columns of empty squares; the $i$-th square under the top square
in the $j$-th
column corresponds to the basic element $T^ix_j$ where $x_j$ is a generator
of the $j$-th direct summand of $V$. The two examples are for the
partitions $(2)$ and $(3,2)$.
$$ V_1=k[T]/(T^2):\quad\beginpicture\setcoordinatesystem units <3mm,3mm>
\multiput {} at 0 -1  1 1 /
\multiput{\sq} at 0 0  0 1 /
\endpicture\qquad
V_2=k[T]/(T^3)\oplus k[T]/(T^2):\quad
\beginpicture\setcoordinatesystem units <3mm,3mm>
\multiput {} at 0 -1  2 2 /
\multiput{\sq} at 0 0  0 1 0 2  1 0  1 1  /
\endpicture
$$

Suppose the quiver $Q$ has {\bf linear orientation,} then a representation
is a triple $(V,U,W)$ where the {\it ambient space} $V$ is
a nilpotent linear operator, the {\it intermediate subspace} $U$ of $V$
is invariant under the action of $T$, and the {\it small subspace}
$W$ of $V$ is also $T$-invariant and contained in the intermediate subspace
$U$.   In each example in this paper, the generators of the intermediate
space can be chosen as basic elements or sums of basic elements of the
ambient space and will be pictured as circles or connected circles that
are contained in the squares for the corresponding basic elements.
Similarly, the generators for the small subspace are indicated by bullets
or connected bullets.
%In the  examples we also list the partition types
%of the ambient, intermediate and small spaces.
$$(V_1, kT, 0): \; \beginpicture\setcoordinatesystem units <4mm,4mm>
\multiput {} at 0 -1  1 1 /
\multiput{\sq} at 0 0  0 1 /
\put{\cir} at 0.5 0
\endpicture %\quad{\rm partitions:} (2), (1), ()
\qquad
(V_2, (T)/(T^3)\oplus (T)/(T^2), k(T^2,T)): \;
\beginpicture\setcoordinatesystem units <4mm,4mm>
\multiput {} at 0 -1  2 2 /
\multiput{\sq} at 0 0  0 1 0 2  1 0  1 1  /
\multiput{\cir} at 0.5 1  1.5 0 /
\multiput{\bul} at 0.5 0  1.5 0 /
\plot 0.5 0  1.5 0 /
\endpicture% \quad{\rm part:} (3,2), (2,1), (1)
$$

\begin{figure}[ht]
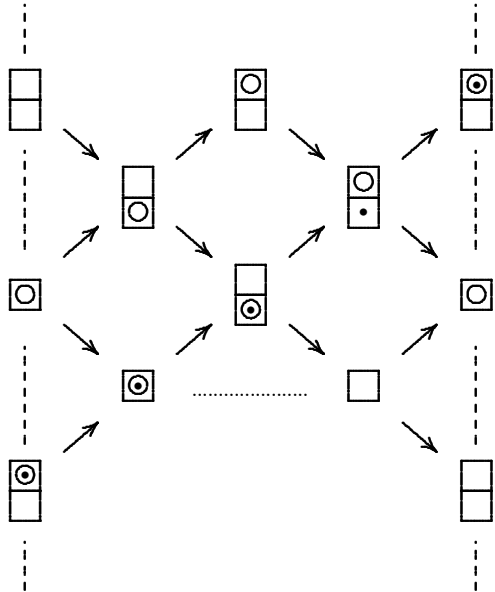

\def\scale#1{\hbox{\beginpicture\setcoordinatesystem units <4mm,4mm>#1
    \multiput{} at -2 -2  2 5 /
  \endpicture}}
$$\arqkLEps$$
\caption{Gorenstein-projective $kQ[\varepsilon]$-modules (linear orientation)}
\end{figure}

For the quiver $Q$ in {\bf V-orientation,} a representation is a triple
$(V,U_\blacktriangleleft,U_\blacktriangleright)$ consisting of the
ambient space $V$ together with a {\it left subspace} $U_\blacktriangleleft$
and a {\it right subspace} $U_\blacktriangleright$, both subspaces are $T$-invariant.
Here is the example which occurs in the middle of the Auslander-Reiten
quiver below.
$$(k[T]/(T^2)\oplus k, k(0,1), k(T,1)): \quad
\beginpicture\setcoordinatesystem units <4mm,4mm>
\multiput {} at 0 -1  2 1 /
\multiput{\sq} at 0 0  0 1  1 0  /
\multiput{\trr} at 1.7 -0.2 /
\multiput{\trl} at 0.3 0.2  1.3 0.2 /
\plot 0.3 0.2  1.3 0.2 /
\endpicture $$

\begin{figure}[ht]
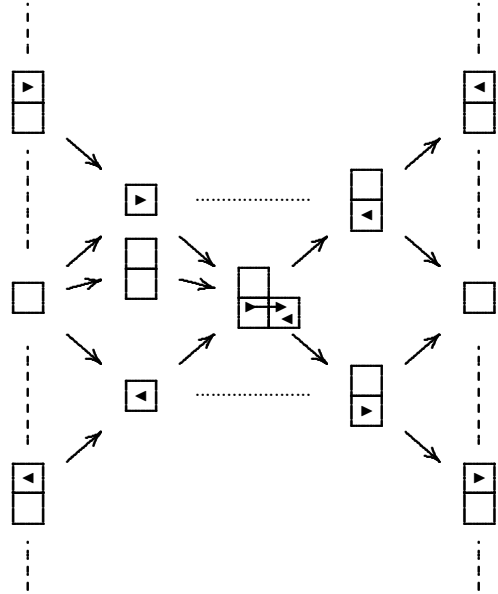

\def\scale#1{\hbox{\beginpicture\setcoordinatesystem units <4mm,4mm>#1
    \multiput{} at -2 -2  2 5 /
  \endpicture}}
$$\arqkVEps$$
\caption{Gorenstein-projective $kQ[\varepsilon]$-modules (V-orientation)}
\end{figure}

For the quiver $Q$ in {\bf $\Lambda$-orientation,} a representation
$(V,W,U_V,U_W)$ consists of two nilpotent linear operators,
the {\it left space} $V$ and the {\it right space} $W$, and also
a {\it subspace} $U_V\cong U_W$
which is embedded as $U_V$ in the left space and as $U_W$ in the right space.
In the
icon, we separate the left space from the right space by a vertical line
and indicate the subspace (abusing notation)
as a subspace of the direct sum of the left
space and the right space.  The example occurs on the left end of the
upper dotted line in the Auslander-Reiten quiver.
$$(k,k[T]/(T^2),k, kT): \quad
\beginpicture\setcoordinatesystem units <3mm,3mm>
\multiput {} at 0 -1  2 1 /
\multiput{\sq} at 0 0  1 0  1 1  /
\multiput{\bul} at 0.5 0  1.5 0 /
\plot 0.5 0  1.5 0 /
\plot 1 -1  1 2 /
\endpicture $$

\begin{figure}[ht]
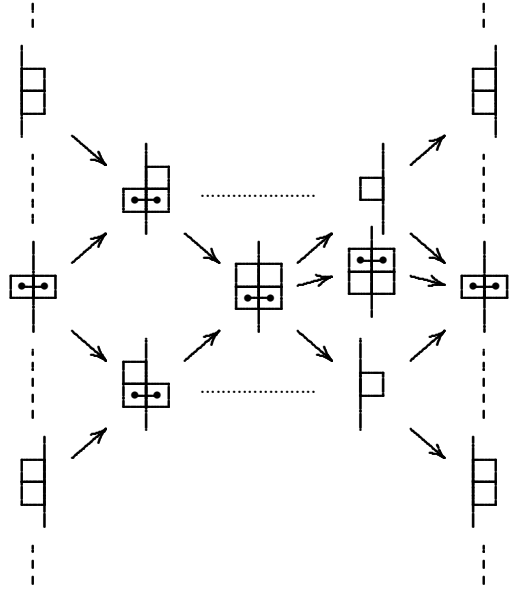

\def\scale#1{\hbox{\beginpicture\setcoordinatesystem units <3mm,3mm>#1
    \multiput{} at -2 -2  2 5 /
    \endpicture}}
$$\arqkLambdaEps$$
\caption{Gorenstein-projective $kQ[\varepsilon]$-modules ($\Lambda$-orientation)}
\end{figure}

%-----------------------------------------------------------------------
  \subsection{Twelve reflections, starting with a simple object}
%-----------------------------------------------------------------------
  \label{sec-12}

  We illustrate the process from Section~\ref{section-reflections}
  which leads to Theorem~\ref{theorem-intro2} in the example where
  $Q$ is the linear quiver of type $\mathbb A_3$,
  $Q=(1\stackrel\alpha\to 2\stackrel\beta\to 3)$,
  and $\Lambda$ any local selfinjective algebra
  with unique maximal ideal $m$ and radical factor
  field $k$.

  \medskip Take $X^0:(0\to 0\to k)$, and reflect at vertex 3:
  In the first step, we compute the cokernel as $\Cok(X^0_\beta)=\bar X^0=k$; in the second step, we take the
  kernel of the projective cover, $Y'_3=m$, and compute $Y'_{\beta^\op}:m\to 0$;
  since this map is not monic, we take in the third step $Y_2=\Lambda$ and $Y_{\beta^\op}:m\to \Lambda$
  the inclusion.  This yields the representation $X^1=(0\to \Lambda\from m)$, for
  $Q(3)=(1\to 2\from 3)$.

  \medskip
  We reflect $X^1$ at vertex 2:  The cokernel of the map
  $(X^1_{\alpha},X^1_{\beta}):X^1_{1}\oplus X^1_{3}\to X^1_{2}$ is $k$, so we obtain $Y'=(0\from m\to m)$,
  where the map $Y'_\beta$ is the identity map, and compute  $X^2=(\Lambda\from m\to m)$.

  \medskip
  We reflect $X^2$ at vertex 1 to obtain  $X^3=(m\to m\to m)$.

  \medskip
  Reflecting $X^3$ at vertex 3 yields $X^4=(m\to m\from 0)$.

  \medskip
  Similarly, we obtain $X^5=(m\from 0\to 0)$; $X^6=(k\to \Lambda\to \Lambda)$,
  $X^7=(k\to \Lambda\from 0)$, $X^8=(k\from k\to\Lambda)$, $X^9=(0\to k\to \Lambda)$,
  $X^{10}=(0\to k\from k)$, $X^{11}=(0\from 0\to k)$, and finally, $X^{12}=X^0=(0\to 0\to k)$,
  so after 12 steps we are back where we started.

%--------------------------------------------------------------------------------
  \subsection{Examples where $T$ has nilpotency index three.}
%--------------------------------------------------------------------------------

  \begin{figure}[ht]
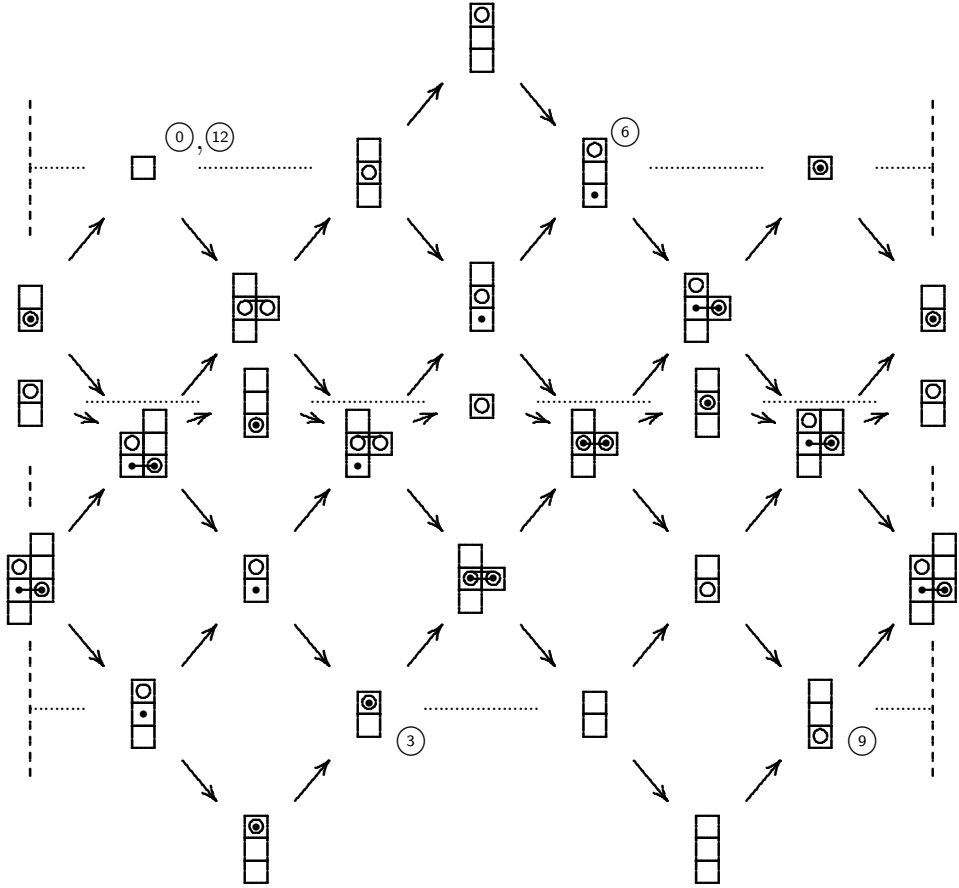

\def\scale#1{\hbox{\beginpicture\setcoordinatesystem units <3mm,3mm>#1
    \multiput{} at 0 -4  1 3 /
  \endpicture}}
\arqkLZeta
\caption{Gorenstein-projective $kQ[\zeta]$-modules (linear orientation)}
\label{figure-zeta-linear}
\end{figure}

  \begin{figure}[ht]
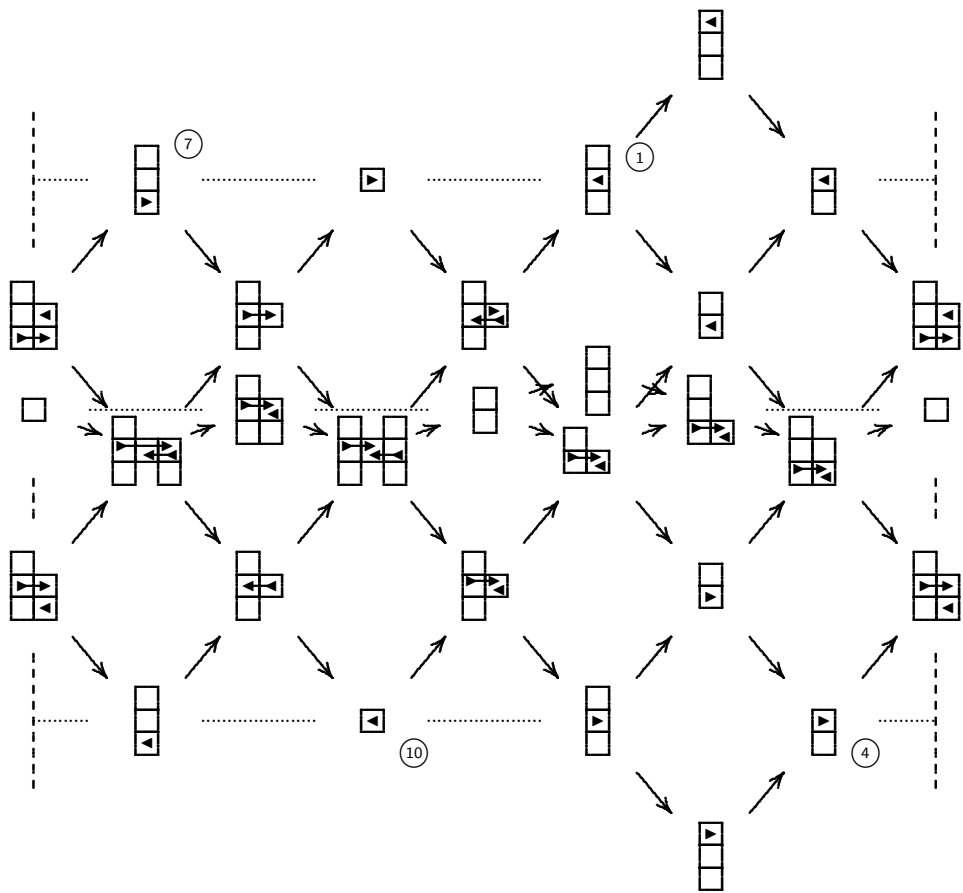

\def\scale#1{\hbox{\beginpicture\setcoordinatesystem units <3mm,3mm>#1
    \multiput{} at 0 -4  3 3 /
  \endpicture}}
\arqkVZeta
\caption{Gorenstein-projective $kQ[\zeta]$-modules (V-orientation)}
\label{figure-zeta-V}
\end{figure}

  \begin{figure}[ht]
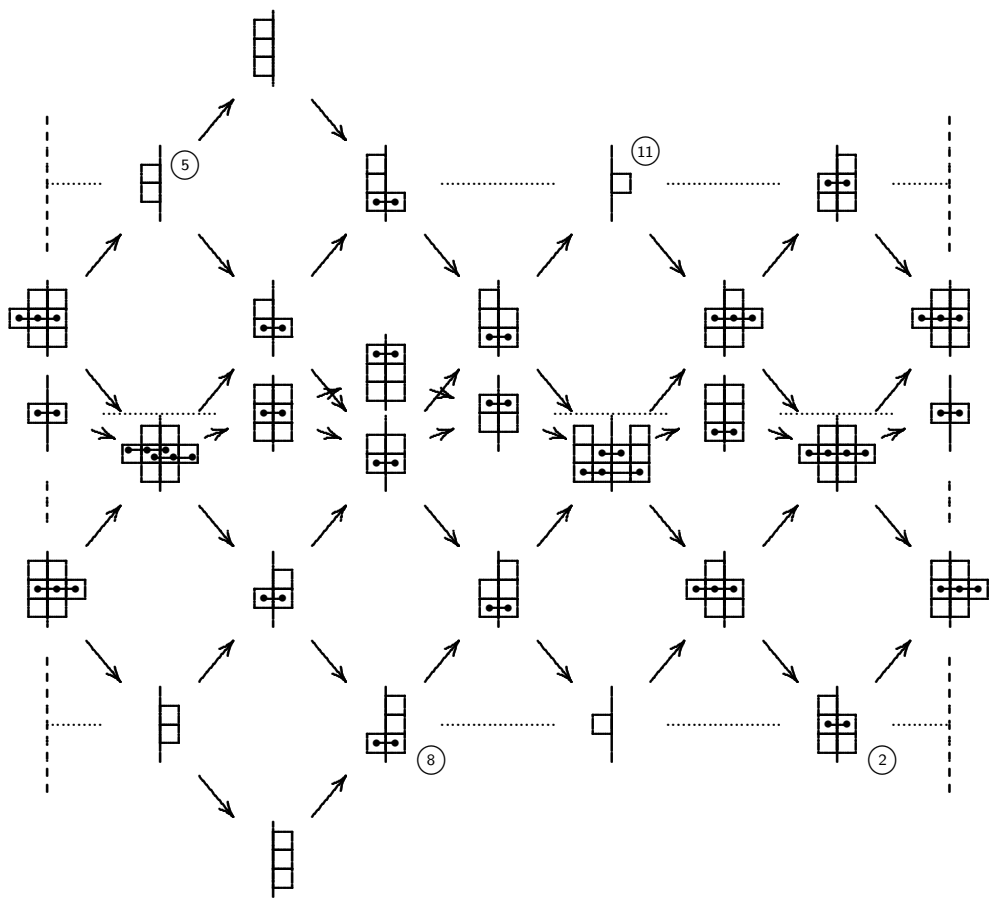

\def\scale#1{\hbox{\beginpicture\setcoordinatesystem units <2.5mm,2.5mm>#1
    \multiput{} at -3 -4 3 3 /
  \endpicture}}
\arqkLambdaZeta
\caption{Gorenstein-projective $kQ[\zeta]$-modules ($\Lambda$-orientation)}
\label{figure-zeta-lambda}
\end{figure}

    Finally, we consider the case where
  $\Lambda=k[\zeta]=k[T]/(T^3)$ and $Q$  is a quiver
  of type $\mathbb A_3$ in either
  linear orientation, V-orientation or $\Lambda$-orientation.

  \medskip
 
    For each of the three orientations, we picture
    in Figures~\ref{figure-zeta-linear}--\ref{figure-zeta-lambda}
    the Auslander-Reiten quiver $\Gamma$
    for $\Gproj\Lambda Q$.  Each has 27 indecomposable Gorenstein-projective modules of which
    24 are stable and 3 are projective-injective.  The shape of the stable part of $\Gamma$ is
    $\overline\Gamma=\mathbb Z\mathbb E_6/(\tau^4)$.  Note that the positions of the three
    projective-injective objects depend on the orientation of $Q$.

  \medskip
  We label each of the above modules $X^i$ in Section~\ref{sec-12} by the symbol \circled i
  in the Auslander-Reiten quiver for the category $\Gproj\Lambda Q$
  where the quiver $Q$
  is oriented such that the module $X^i$ is defined.

    \medskip
    We have seen that each reflection gives rise to a categorical equivalence
    of the stable category of Gorenstein-projective modules, hence gives rise to an isomorphism
    of translation quivers between the stable parts of the AR-quivers.
    It turns out that in this case, this isomorphism is determined
    uniquely by each pair
    $X^i$, $X^{i+1}$ of successive objects.
    Thus we can read off from the AR-quivers
    how the reflection functors act on each object.

    \medskip
    In particular the orbit which is pictured in Figure~\ref{figure-one-orbit}
    is obtained by tracing the objects on the central $\tau$-orbits in Figures~\ref{figure-zeta-linear}--\ref{figure-zeta-lambda}
    under the action of those translation quiver isomorphisms.

\bigskip
%\newpage
Address of the authors:

\parbox[t]{5.5cm}{\footnotesize\begin{center}
              Department of Mathematics\\
              Nantong University\\
              Nantong 226019\\
              China\end{center}}
\parbox[t]{5.5cm}{\footnotesize\begin{center}
              Mathematical Sciences\\
              Florida Atlantic University\\
              Boca Raton, Florida 33431\\
              United States of America\end{center}}

\smallskip \parbox[t]{5.5cm}{\centerline{\footnotesize\tt xiuhualuo@ntu.edu.cn}}
           \parbox[t]{5.5cm}{\centerline{\footnotesize\tt markus@math.fau.edu}}

\end{document}